\documentclass[letterpaper,11pt]{article}

\title{Algebraic Geometry of Bayesian Networks}

\author{Luis David Garcia, Michael Stillman and Bernd Sturmfels}
\date{}

\usepackage{latexsym,array,delarray,amsthm,amssymb}
\usepackage[sumlimits]{amsmath}
\usepackage[dvips]{graphicx}

\theoremstyle{plain}
\newtheorem{thm}{Theorem}

\newtheorem{prop}[thm]{Proposition}
\newtheorem{cor}[thm]{Corollary}
\newtheorem{conj}[thm]{Conjecture}

\theoremstyle{definition}

\newtheorem{ex}[thm]{Example}

\theoremstyle{remark}

\newcommand{\zz}{\mathbb{Z}}

\newcommand{\pp}{\mathbb{P}}

\newcommand{\rr}{\mathbb{R}}
\newcommand{\cc}{\mathbb{C}}

\newcommand{\mm}{\mathcal{M}}

\newcommand{\rto}{\longrightarrow}
\newcommand{\lto}{\longleftarrow}

\newcommand{\perpp}{\mbox{$\perp\!\!\!\perp$}}

\DeclareMathOperator{\local}{local}
\DeclareMathOperator{\pa}{pa}
\DeclareMathOperator{\nd}{nd}
\DeclareMathOperator{\gl}{global}

\begin{document}
\maketitle

\begin{abstract}
\noindent
We study the algebraic varieties defined by the
conditional independence statements of Bayesian
networks. A complete algebraic classification is given for
Bayesian networks on at most five random variables.
Hidden variables are related to the
geometry of  higher secant varieties.
\end{abstract}

\section{Introduction}

The emerging field of {\em algebraic statistics}  \cite{PRW}
advocates polynomial algebra as a tool
in the statistical analysis of
experiments and discrete data.
Statistics textbooks define a {\em statistical model}
as a family of probability distributions, and a closer look
reveals that these families are often
 real algebraic varieties: they are the
zeros of some polynomials in the probability simplex
\cite{GHKM}, \cite{SS1}.

In this paper we examine {\em directed graphical models} for
discrete random variables. Such models are also known
as {\em Bayesian networks} and they are widely used
in machine learning, bioinformatics and many other
applications \cite{L}, \cite{P}. Our aim is to place
Bayesian networks into the realm of algebraic statistics, by
developing the necessary theory in algebraic geometry
and by demonstrating the effectiveness of
Gr\"obner bases for this class of models.

Bayesian networks can be described in
two possible ways, either by a recursive factorization of
probability distributions or by conditional independence
statements (local and global Markov properties).
This is an instance of the computer algebra principle
that varieties can be presented
either parametrically or implicitly \cite[\S 3.3]{CLO}.
The equivalence of these two representations
for Bayesian networks is a well-known theorem
in  statistics \cite[Theorem 3.27]{L},
but, as we  shall see,  this theorem
is surprisingly delicate and no longer holds when
probabilities are replaced by negative reals
or complex numbers. Hence in the usual setting of algebraic
geometry, where the zeros lie in $\cc^d$, there are many
``distributions'' which satisfy the global Markov
property but which do not permit a recursive  factorization.
We explain this phenomenon using primary  decomposition of
polynomial ideals.

This paper is organized as follows.
In Section 2 we review the algebraic
theory of conditional independence, and we
explicitly determine the  Gr\"obner basis and primary
decomposition arising from the {\em contraction axiom}
\cite{P}, \cite[\S 2.2.2]{Stu}.
This axiom is shown to fail for negative real numbers.
In Section 3 we introduce the ideals $I_{{\rm local}(G)}$ and
$I_{{\rm global}(G)}$ which represent a
a Bayesian network $G$. When $G$ is a forest
then these ideals are the toric ideals
derived from  undirected graphs as in \cite{GMS};
see Theorem \ref{treethm} below.

The recursive factorization of a Bayesian network gives
rise to a map between polynomial rings which
is studied  in Section 4.
The kernel of this {\em factorization map} is the
{\em distinguished prime ideal}. We prove that this prime
is always a reduced primary component of
$I_{{\rm local}(G)}$ and $I_{{\rm global}(G)}$.
Our results in that section include the
solutions to Problems 8.11 and 8.12 in  \cite{S}.

In Sections 5 and 6 we present the results of our
computational efforts: the complete algebraic classification
of all Bayesian networks on four arbitrary random variables
and all Bayesian networks on  five binary random variables.
The latter involved computing the primary
decomposition of $301$ ideals generated by
a large number of quadrics in $32$ unknowns.
These large-scale primary decompositions
were carried out in {\tt Macaulay2} \cite{GS}. Some of the
techniques and software tools we used are described in the Appendix.

The appearance of hidden variables in
 Bayesian networks leads to challenging problems
in  algebraic geometry. Statisticians have known for decades  that
the dimension of the corresponding varieties
can unexpectedly drop \cite{Good},
but the responsible singularities have been
studied only quite recently, in
\cite{GHKM} and \cite{SS1}. In Section 7
we examine the elimination problem arising
from hidden random variables, and we
relate it to problems in  projective algebraic geometry.
We demonstrate that the {\em naive Bayes model}
corresponds to the higher secant varieties
of Segre varieties (\cite{CJ}, \cite{CGG}),
and we present several new results on
the dimension and defining ideals of these secant
varieties.

Our algebraic theory does not compete with
but rather complements other approaches to
conditional independence models. An
 impressive combinatorial theory of such models
 has been developed by
Mat\'u\v s \cite{Mat} and Studen\'y \cite{Stu},
culminating  in their characterization of all
realizable independence models on four
random variables. Sharing many of the views expressed by these
authors, we believe that exploring the
precise relation between their work and ours
will be a very fruitful research direction for the near future.

\bigskip

\noindent {\bf Acknowledgements.}
Garcia and Sturmfels were partially supported
by the CARGO program of the National Science Foundation
(DMS-0138323). Stillman was partially supported
by  NSF grant DMS 9979348,
and  Sturmfels was also partially supported by
NSF grant DMS-0200729.

\section{Ideals, Varieties and Independence Models}

We begin by reviewing the general algebraic framework
for independence models presented in \cite[\S 8]{S}.
Let $X_{1}, \dots , X_{n}$ be discrete random variables
where $X_{i}$ takes values in the
finite set $[d_{i}] = \{1,2,\dots,d_{i}\}$. We
write $D = [d_{1}]\times [d_{2}]\times \dotsm \times [d_{n}]$ so that
$\rr^{D}$
denotes the real vector space of $n$-dimensional tables of format
$d_{1}\times \dotsm \times d_{n}$. We introduce an indeterminate
$p_{u_{1}u_{2}\dotsm u_{n}}$ which represents the probability of the
event $X_{1} = u_{1},\, X_{2} = u_{2},\dots ,X_{n} = u_{n}$. These
indeterminates generate the ring $\rr[D]$ of polynomial functions on
the space of tables $\rr^{D}$.
A \emph{conditional independence statement}
has the form
\begin{equation}\label{ind-state}
A \textrm{ is independent of } B \textrm{ given } C \qquad (\textrm{
   in symbols: } A \perpp B \mid C)
\end{equation}
where $A, B$ and $C$ are pairwise disjoint subsets of $\{X_{1},\dots
,X_{n}\}$. If $C$ is empty then \eqref{ind-state} means that
$A$ is independent of $B$.
By \cite[Proposition 8.1]{S},
the  statement \eqref{ind-state} translates into a set
   of homogeneous quadratic polynomials in $\rr[D]$,
and we write $I_{A \perpp B \mid C}$ for the
ideal generated by these polynomials.

Many  statistical models (see e.g.~\cite{L},  \cite{Stu})
  can be described by a
finite set of independence statements \eqref{ind-state}. An
\emph{independence model} is any such set:
\[\mm \quad = \quad
  \{A^{(1)} \perpp B^{(1)} \mid C^{(1)}, \dots, A^{(m)} \perpp B^{(m)}
\mid C^{(m)}\}. \]
The ideal of the independence model $\mm$ is defined
as the sum of ideals
\[I_{\mm} \quad = \quad I_{A^{(1)} \perpp B^{(1)}
\mid C^{(1)}} + \dotsm + I_{A^{(m)} \perpp B^{(m)} \mid C^{(m)}} .  \]
We wrote code in {\tt Macaulay2} \cite{GS} and {\tt Singular} \cite{GPS}
 for generating the ideals $I_{\mm}$.
The {\em independence variety} is the
set $V(I_{\mm})$ of common zeros in $\cc^D$
of the polynomials in $I_\mm$. Equivalently,
$V(I_{\mm})$ is the set of all
  $d_{1}\times \dotsm \times d_{n}$-tables with complex
number entries which satisfy the
conditional independence statements in $\mm$.
The variety $V(I_{\mm})$ has three natural subsets:
\begin{itemize}
\item the subset of real tables, denoted $\, V_{\rr}(I_\mm)$,
\item the non-negative tables, denoted $\,  V_{\geq}(I_\mm)$,
\item the non-negative tables whose entries
sum to one, $\, V_{\geq}(I_\mm + \langle p-1\rangle)$,
\end{itemize}
Here $\,p\,$ denotes the sum of all unknowns
$p_{u_1 \cdots u_n}$, so that
$\, V_{\geq}(I_\mm + \langle p-1\rangle)\,$
is the subset of the probability simplex
specified by the model $\mm$.

We illustrate these definitions by
analyzing the independence model
$\,\, {\cal M}  \,\, = \,\, \bigl\{
\,\, 1 \perpp 2 \,\, | \,\, 3 \, , \,\,
\,\, 2 \perpp 3 \,\bigr\}\,$ for $n=3$
discrete random variables. Theorem \ref{threemodel}
will be cited in Section 5 and it serves
as a preview to Theorem \ref{binaryfour}.
 The ideal $I_{\cal M}$ lies in the polynomial
ring $\rr[D]$ in  $\,d_1 d_2 d_3 \,$ unknowns $\, p_{i j k}$.
Its minimal generators are
$\, \binom{d_1}{2}\binom{d_2}{2} d_3\,$ quadrics of the form
$\,p_{i j k}\, p_{r s k} -  p_{i s k} \, p_{r j k}  \,$ and
$\, \binom{d_2}{2}\binom{d_3}{2}\,$ quadrics
of the form $\,p_{+ j k} \, p_{+ s t} -  p_{+ j t} \, p_{+ s k} $.
We change coordinates in $\rr[D]$
by replacing each unknown $\,p_{1 j k} \,$ by
$\, p_{+ j k } \, = \,  \sum_{i=1}^{d_1} p_{i j k} $.
This coordinate change transforms $I_{\cal M}$
into a binomial ideal in $\rr[D]$.

\begin{thm}
\label{threemodel}
The ideal $I_{\cal M}$ has a Gr\"obner basis consisting of
squarefree binomials of degree two, three and four, and
it is hence radical. It has $2^{d_3} \! - \! 1$ minimal primes,
each generated by the $2 \times 2$-minors
of a generic matrix.
\end{thm}

\begin{proof}
The minimal primes of $I_{\cal M}$ will be indexed by
proper subsets of $\,[d_3]$. For each such subset
$\sigma$ we introduce the monomial prime
$$ M_\sigma \quad = \quad \langle
\, p_{+jk} \,\,\, | \,\,\,
  \, j \in [d_2], \, k \in \sigma \,\rangle , $$
and the complementary monomial
$$ m_\sigma  \quad = \quad
\prod_{j=1}^{d_2}  \prod_{k \in [d_3] \backslash \sigma} p_{+jk}\,, $$
and we define the ideal
$$ P_\sigma \quad := \quad
\bigl((I_{\cal M} \,+\, M_\sigma) : m_\sigma^\infty \bigr). $$
It follows from the general theory of binomial ideals
\cite{ES} that $P_\sigma$ is a binomial prime ideal.
A closer look reveals that $P_\sigma$ is minimally generated by the
$\,d_2 \cdot |\sigma| \,$ variables in $M_\sigma$
together with all the $2 \times 2$-minors of the following
two-dimensional matrices: the matrix $\,(\, p_{ijk}\,)\,$ where the
rows are indexed by $j \in [d_2]$ and the columns are indexed by pairs
$(i,k)$ with
$\,i \in \{+,2,3,\ldots,d_1\}$ and $\, k \in [d_3] \backslash
\sigma$, and for each $k \in \sigma$, the matrices $\,(\, p_{ijk}\,)\,$
where the
rows are indexed by $j \in [d_2]$ and the columns are indexed by
$\,i \in \{2,3,\ldots,d_1\}$.

We partition $\, V(I_{\cal M})\,$
into $2^{d_3}$ strata, each indexed by a subset $\sigma$
of $[d_3]$. Namely, given a point $(p_{ijk})$ in $V(I_{\cal M})$ we
define the subset $\sigma$ of $[d_{3}]$ as the set
of all indices $k$ such that
$\,(p_{+1k}, p_{+2k}, \ldots,p_{+d_2 k})\,$
is the zero vector. Note that two tables $(p_{ijk})$ lie in the same
stratum if and only if they give the same $\sigma$.
The stratum indexed by $\sigma$ is a dense subset
in  $\,V(P_\sigma) $.
When $\sigma = [d_{3}]$ the stratum consists of all tables such that
the line sums $p_{+jk}$ are all zero, and for each fixed $k$, the
remaining $(d_{1}-1)\times d_{2}$-matrix $(p_{ijk})$ with $i\geq 2$
has rank $\leq 1$. So this locus is defined by the prime ideal
$P_{[d_{3}]}$. Any point in this stratum
satisfies the defining equations of $P_{\sigma}$ for any proper subset
$\sigma$. So
the stratum indexed by $[d_3]$ lies in the closure
of all other strata. But all remaining $2^{d_3}-1$ strata have
the property that no stratum lies in the
closure of any other stratum, since
the generic point of $P_{\sigma}$ lies in exactly one stratum for any
proper subset $\sigma$. Hence $\,V(I_{\cal M})\,$ is the
irredundant union of the irreducible varieties
$\, V(P_\sigma)\,$ where $\sigma$ runs
over all proper subsets of $[d_3]$.
The second assertion  in Theorem \ref{threemodel}
now follows from Hilbert's Nullstellensatz.

To prove the  first assertion, let us first
note that $\,P_\emptyset \,$ is the
prime ideal of $2 \times 2 $-minors
of the $\,d_2 \times (d_1 d_3)$-matrix
$\,(p_{ijk})\,$ with rows  indexed
by $j \in [d_2]$ and columns
indexed by pairs $(i,k) \in \{+,2,3,\ldots,d_1\} \times [d_3]$. Hence
\begin{equation}
\label{impliii}
  P_\emptyset \quad = \quad
\bigl(\,I_{\cal M} \,: \, m_{\emptyset}^\infty \,\bigr)
\quad = \quad
I_{2 \,\perpp\, \{1,3\}}.
\end{equation}
It is well known (see e.g.~\cite[Proposition 5.4]{GBCP})
that the quadratic generators
\begin{equation}
\label{knowngb}
p_{ijk}\, p_{rst} \quad - \quad p_{isk}\, p_{rjt}
\end{equation}
form a reduced Gr\"obner basis for (\ref{impliii}) with respect
to the ``diagonal term order''. We modify this
Gr\"obner basis to a Gr\"obner basis for $I_{\cal M}$
as follows:

\begin{itemize}
\item if $\,k=t \, $ \ take  \ (\ref{knowngb}),
\item if $\, i = + \,$ and $\, r = + \,$ \  take \  (\ref{knowngb}),
\item if $\,i = + \,$ and $\,r \not= +\,$
and $k \not= t$ \ take \ (\ref{knowngb}) times
$p_{+jt}$\  for any $j$,
\item if $\,i \not= + \,$ and $\,r \not= +\,$
and $k \not= t$ \ take \ (\ref{knowngb}) times
$p_{+jt} p_{+sk}$\  for any $j,s$.
\end{itemize}
All of these binomials lie in $I_{\cal M}\,$
(this can be seen by taking S-pairs of the generators)
and their S-pairs reduce to zero. By Buchberger's
criterion, the given set of quadrics, cubics and quartics is
a Gr\"obner basis, and the corresponding initial monomial
ideal is square-free.
This implies that $I_{\cal M}$
is radical (by \cite[Proposition 5.3]{S}), and the proof  is complete.
\end{proof}

The theorem above can be regarded as an algebraic refinement
of the following well-known rule
for conditional independence
(\cite{P}, \cite[\S 2.2.2]{Stu}).

\begin{cor} {\rm (Contraction Axiom) } If a probability
distribution on $[d_1]\times [d_2] \times [d_3]$
satisfies $\,1 \perpp 2 \,|\, 3 \,$
and $\,2 \perpp 3\,$ then it also satisfies
$\,2 \perpp \{1,3\}$.
\end{cor}

\begin{proof}
The non-negative points satisfy
$\, V_{\geq}(P_\sigma) \subseteq V_{\geq}(P_\emptyset)$,
and this implies
$$ V_{\geq}(I_{\cal M}) \quad = \quad
V_{\geq}(I_{2 \,\perpp\, \{1,3\}} ). $$
Intersecting with the probability simplex yields the assertion.
\end{proof}

Theorem  \ref{threemodel} shows
that the Contraction Axiom fails to hold when
probabilities are replaced by negative real numbers.
Any general point on  $V(P_\sigma)$ for $\sigma \not= \emptyset \,$
satisfies $\,1 \perpp 2 \,|\, 3 \,$
and $\,2 \perpp 3\,$ but it does not satisfy
$\,2 \perpp \{1,3\}$.

\section{Algebraic Representation of Bayesian Networks}

A \emph{Bayesian network} is an acyclic directed graph
$G$ with vertices $X_{1},\dots,X_{n}$.   The following notation
and terminology is consistent with
Lauritzen's book \cite{L}.
The  \emph{local Markov property}
on $G$ is the set of independence
statements
\[\local(G) \quad = \quad \{X_{i} \perpp \nd(X_{i}) \mid \pa(X_{i}) \:
\colon  \:
i = 1,2,\dots,n\}, \]
where $\pa(X_{i})$ denotes the
set of parents of $X_{i}$ in $G$ and $\nd(X_{i})$ denotes the
set of nondescendents of $X_{i}$ in $G$. Here $X_j$ is
a \emph{nondescendent} of $X_i$ if there is no directed
path from $X_i$ to $X_j$ in $G$.
The \emph{global Markov
property}, $\gl(G)$, is the
set of independence statements $A \perpp B \mid C$, for any triple
$A,B,C$ of subsets of vertices of $G$ such that $A$ and $B$
are \emph{d-separated} by $C$.
Here two subsets $A$ and $B$ are said to be \emph{d-separated} by $C$ if
all chains from $A$ to $B$ are blocked by $C$.
A chain $\pi$ from $X_{i}$ to $X_{j}$ in $G$
is said to be \emph{blocked} by
a set $C$ of nodes if it contains a vertex $X_{b}
\in \pi$ such that either
\begin{itemize}
   \item $X_{b} \in C$ and arrows of $\pi$ do not meet head-to-head at
   $X_{b}$, or
\item $X_{b} \notin C$ and $X_{b}$ has no descendents in $C$, and
   arrows of $\pi$ do meet head-to-head at $X_{b}$.
\end{itemize}

For any Bayesian network $G$, we have
$\, \local(G) \, \subseteq \, \gl(G) $, and this
implies the following containment relations between ideals
and varieties
\begin{equation}\label{eq:ideal}
I_{\local(G)} \, \subseteq \, I_{\gl(G)} \qquad \hbox{and} \qquad
V_{\local(G)} \, \supseteq \, V_{\gl(G)}.
\end{equation}
The latter inclusion extends to the three real varieties listed
above, and we shall discuss when equality holds.
First, however, we give an algebraic version of the
description of Bayesian networks by recursive factorizations.

Consider the set of parents
of the $j$-th node,
$\pa(X_{j}) = \{X_{i_{1}},\dots,
X_{i_{r}}\}$,
  and consider any event 
$X_{j} = u_{0}$ conditioned on
$X_{i_{1}} = u_{1}, \dots, X_{i_{r}} = u_{r}$,
where $1 \leq u_{0} \leq d_{j}, 1 \leq u_{1} \leq d_{i_{1}}, \dots, 1
\leq
u_{r} \leq d_{i_{r}}$. We introduce an unknown
$\,q^{(j)}_{u_{0}u_{1}\dotsm u_{r}}$ to denote the
conditional probability of
this event, and we subject these unknowns to the linear relations
$\,\sum_{v=1}^{d_{j}} q^{(j)}_{vu_{1}\dotsm u_{r}} = 1 \,$
for all  $ 1 \leq u_{1} \leq d_{i_{1}}, \dots, 1 \leq
u_{r} \leq d_{i_{r}}$. Thus, we have introduced
$(d_{j}-1)d_{i_{1}}\dotsm d_{i_{r}}$ unknowns for the vertex $j$.
Let $E$ denote the set of these unknowns
$q^{(j)}_{u_{0}u_{1}\dotsm  u_{r}}$ for all $j \in \{1,\ldots,n\}$,
and let $\rr[E]$ denote the polynomial ring they generate.

If the $n$  random variables are binary
($d_i = 2$ for all $i$) then the notation for $\rr[E]$
can be simplified by dropping the first lower index
and writing
\[q^{j}_{u_{1}\dotsm u_{r}}
\quad := \quad q^{(j)}_{1u_{1}\dotsm u_{r}}
\quad = \quad 1 \, - \, q^{(j)}_{2u_{1}\dotsm u_{r}}\]
In the binary case, $\rr[E]$ is a polynomial ring in
$\,\sum_{j=1}^{n}2^{|\pa(X_{j})|} \,$ unknowns.

The factorization of  probability distributions according to $G$
defines a polynomial map
$\,\phi \, : \, \rr^E \rightarrow \rr^D$.
By restricting to non-negative reals we get an
induced map $\phi_{\geq 0}$.
These maps are specified by the ring homomorphism
$ \Phi \,: \, \rr[D] \rightarrow \rr[E]\,$ which
takes the unknown
$\,p_{u_1 u_2 \cdots u_n} \,$ to the product
of the expressions
$\,q^{(j)}_{u_j u_{i_1}\dotsm u_{i_r}}  \,$
as $j$ runs over $\{1,\ldots,n\}$.
The image
of $\phi$ lies in the independence variety
$\,V_{{\rm global}(G)} $, or, equivalently, the independence ideal
$I_{{\rm global}(G)}$ is contained in the prime ideal $\, \ker(\Phi)$.
The \emph{Factorization Theorem} for Bayesian networks
\cite[Theorem 3.27]{L} states:

\begin{thm}
\label{factheorem}
The following four subsets of the probability simplex coincide:
\begin{eqnarray*}
& V_{\geq}(I_{\local(G)} + \langle p -1 \rangle) \,\,\, =\,\,\,
  V_{\geq}(I_{\gl(G)} + \langle p - 1 \rangle) \\
& =\quad   V_{\geq}(\ker(\Phi))\quad = \quad
{\rm image}(\phi_{\geq}).
\end{eqnarray*}
\end{thm}

\begin{ex}
Let $G$ be the network on
three binary random variables which
has a single directed
edge from $3$ to $2$. The parents and nondescendents are
$$
{\rm pa}(1) \! = \! \emptyset, {\rm nd}(1) \! =\! \{2,3\},\,
{\rm pa}(2) \! = \! \{3 \}, {\rm nd}(2) \!= \! \{1\},\,
{\rm pa}(3) \! = \! \emptyset, {\rm nd}(3) = \{1\}. $$
The resulting conditional independence statements are
$$ {\rm local}(G) \,=\,  {\rm global}(G) \, =\,
\bigl\{ 1 \perpp 3,\; 1\perpp 2\mid 3,\; 1\perpp \{2,3\} \bigr\}. $$
The ideal expressing the first two statements is contained
in the ideal expressing the third statement, and we find
that $\, I_{{\rm local}(G)} \,=\,
I_{ 1\perpp \{2,3\} }\,$ is the
ideal generated by the  six $2 \times 2$-subdeterminants
of the $2\times 4$-matrix
\begin{equation}
\label{twobyfour}
\begin{array}({cccc})
p_{111} & p_{112} & p_{121} & p_{122} \\
p_{211} & p_{212} & p_{221} & p_{222}
\end{array}
\end{equation}
This ideal is prime and its generators form a Gr\"obner basis.
The Factorization Theorem is understood as follows
for this example. We have
$\,E\, = \{  q^{1},  q^{2}_1,  q^{2}_2,  q^{3}\} $,
and our ring map $\Phi$ takes
the matrix (\ref{twobyfour}) to
$$
\begin{array}({cccc})
    \! \!  q^{1} q^{2}_1 q^{3} &     q^{1} q^{2}_2 (1\!-\!q^{3}) &
     q^{1} (1\!-\!q^{2}_1) q^{3} &     q^{1} (1\!-\!q^{2}_2)
(1\!-\!q^{3})
   \!  \! \\
   \!\!
(1\!-\!q^{1}) q^{2}_1 q^{3} & (1\!-\!q^{1}) q^{2}_2 (1\!-\!q^{3}) &
(1\!-\!q^{1}) (1\!-\!q^{2}_1) q^{3} & (1\!-\!q^{1}) (1\!-\!q^{2}_2)
(1\!-\!q^{3})   \! \! \!
\end{array}
$$
The map $\phi$ from $\rr^4$ to $\rr^8$ corresponding to
the ring map $\Phi : \rr[D] \rightarrow \rr[E]$ gives
a parametrization of all $2 \times 4$-matrices of rank $1$
whose entries sum to $1$. The Factorization Theorem for $G$
is the same statement for non-negative matrices. The kernel
of $\Phi$ is exactly equal to $\, I_{{\rm local}(G)} + \langle p -1
\rangle$. \qed
\end{ex}

\vskip .1cm

Our aim is to decide to what extent the
Factorization Theorem is valid over all
real and all complex numbers. The corresponding algebraic question
is to study the ideal $I_{{\rm local}(G)}$ and
to determine its primary decomposition.
Let us begin by considering all Bayesian networks on
three random variables. We shall prove that
for such small networks the ideal $I_{{\rm local}(G)}$
is always prime and coincides with the kernel of $\Phi$.
The following theorem is  valid for
arbitrary positive integers $d_1,d_2,d_3$.
It is not restricted to the binary case.

\begin{prop}
\label{threeallprime}
For any Bayesian network $G$ on three
discrete random variables, the
ideal $I_{{\rm local}(G)}$ is prime, and it
has a quadratic Gr\"obner basis.
\end{prop}

\begin{proof} We completely classify all possible cases.
  If $G$ is the complete graph, directed acyclically,
then  ${\rm local}(G)$ contains no
nontrivial independence statements, so
$I_{{\rm local}(G)}$ is the zero ideal.
In what follows we always exclude this case.
There are five isomorphism types
of (non-complete) directed acyclic graphs
on three nodes. They correspond to the rows
of the following table:

\begin{table}[!htb]
\begin{center}
\begin{tabular}{| l | c | c |}
\hline
Graph & Local/Global Markov property & Independence ideal \\ \hline
$\begin{array}{ccc} 3 & 2 & 1 \end{array}$ &
$1 \perpp \{2,3\},\; 2 \perpp \{1,3\},\; 3 \perpp \{1,2\}$
& $I_{\mathrm{Segre}}$
\\ \hline
$3 \rto 2 \quad 1$ & $1 \perpp 3,\; 1\perpp 2\mid 3,\; 1\perpp \{2,3\}$ &
$I_{1\perpp \{2,3\}}$
\\ \hline
$3 \rto 2 \rto 1$ & $1\perpp 3 \mid 2$ & $I_{1 \perpp 3\mid 2}$
\\ \hline
$1 \lto 3 \rto 2$ & $1\perpp 2 \mid 3$ & $I_{1 \perpp 2\mid 3}$
\\ \hline
$3 \rto 1 \lto 2$ & $2\perpp 3$ & $I_{2\perpp 3}$
\\ \hline
\end{tabular}
\end{center}
\end{table}

The third and fourth network represent the
same independence model.  In all cases except for the first,
the ideal $I_{{\rm local}(G)}$ is of the form
$\,I_{A \perpp B \mid C}$, i.e., it is specified
by a single independence statement. It was shown in
\cite[Lemma 8.2]{S} that such ideals are prime.
They are determinantal ideals and well known
to possess a quadratic Gr\"obner basis.
The only exceptional graph is the empty graph,
which leads to the model of complete independence
$\,1 \perpp \{2,3\},\; 2 \perpp \{1,3\},\; 3 \perpp \{1,2\}$.
The corresponding ideal defines the
Segre embedding of the product of three
projective spaces $\,\pp^{d_1-1} \times \pp^{d_2-1} \times
\pp^{d_3-1}\,$
into $\,\pp^{d_1 d_2 d_3 - 1}$.  This ideal is prime
and has a quadratic Gr\"obner basis.
\end{proof}

A network $G$ is a {\em directed forest} if
every node has at most one parent.
The conclusion of Proposition \ref{threeallprime}
also holds for directed forests
on any number of nodes. Proposition
\ref{wefoundatree} will show that
the direction of the edges is crucial:
it is not sufficient to assume that
the underlying undirected graph is a forest.

\begin{thm}
\label{treethm}
Let $G$ be a directed forest.
Then  $I_{\mathrm{global}(G)}$ is prime and has a
quadratic Gr\"obner basis. These
properties generally fail for $I_{\mathrm{local}(G)}$.
\end{thm}

\begin{proof}
For a direct forest,  the definition of a \emph{blocked} chain
reads as follows. A chain $\pi$ from $X_i$ to $X_j$ in $G$ is blocked
by a
set $C$ if it contains a vertex $X_b \in \pi \cap C$. Hence, $C$
d-separates $A$ from $B$ if and only if  $C$ separates $A$ from $B$ in
the
undirected graph underlying $G$. Thus, \cite[Theorem 12]{GMS} implies
that
$\,I_{{\rm global}(G)}\,$ coincides with the distinguished prime ideal
${\rm ker}(\Phi)$, this ideal has a quadratic Gr\"obner basis.
The second assertion is proved by the networks
{\tt 18} and {\tt 26}
in Table \ref{binary4}. See also \cite[Example 8.8]{S}.
\end{proof}

We close this section with a conjectured
characterization of  the global Markov
property on a Bayesian network $G$
in terms of commutative algebra.

\begin{conj}\label{conj:quadrics}
$I_{{\rm global}(G)}$ is the ideal generated by all quadrics in
$\, \ker(\Phi)$.
\end{conj}

\section{The Distinguished Component}

In what follows we shall assume that every edge
$(i,j)$ of the Bayesian network $G$ satisfies $i > j$.
In particular, the node $1$ is always a sink and
the node $n$ is always a source.
For any integer $r \in [n]$ and $u_i \in [d_i]$ as before,
we abbreviate the \emph{marginalization} over the
first $r$ random variables as follows:
$$ p_{++\cdots+ u_{r+1} \cdots u_n} \quad := \quad
\sum_{i_1=1}^{d_1}
\sum_{i_2=1}^{d_2} \cdots
\sum_{i_r=1}^{d_r}
p_{i_1 i_2 \cdots i_r u_{r+1} \cdots u_n}. $$
This is a linear form in our polynomial ring $\rr[D]$.
We denote by ${\bf p}$ the product of
all of these linear forms. Thus the
equation of $\,{\bf p} = 0\,$ defines
a hyperplane arrangement in $\rr^D$.
We shall prove that the ideal $I_{{\rm local}(G)}$ is prime
locally outside this hyperplane arrangement,
and hence so is $I_{{\rm global}(G)}$.
The following theorem provides the solution to
\cite[Problem 8.12]{S}.

\begin{thm}
\label{prob812solved}
  The prime ideal   $\,\ker(\Phi)\,$
  is a minimal primary component of
both of the ideals
  $I_{\local(G)}$ and $ I_{\gl(G)}$. More precisely,
\begin{equation}
\label{sateqn}
  \bigl( I_{\local(G)} : {\bf p}^\infty \bigr)    \,\,\, = \,\,\,
    \bigl( I_{\gl(G)} : {\bf p}^\infty \bigr)  \,\,\, = \,\,\,
\ker(\Phi).
\end{equation}
\end{thm}

The prime ideal $\, \ker(\Phi) \,$ is called  the
\emph{distinguished component}. It can be characterized
as the set of all
homogeneous polynomial functions on $\rr^D$ which
vanish on all probability distributions that
factor according to $G$.

\begin{proof}
We relabel $G$ so that ${\rm pa}(1) = \{2,3,\ldots,r\}$ and
${\rm nd}(1) = \{r+1,\ldots,n\}$.
Let $A$ denote a set of  $(d_1-1)d_2 \cdots d_r$ new
unknowns $\,a_{i_1 i_2 \cdots i_r} $, for $i_1 > 1$
defining a polynomial ring  $\rr[A]$.  Define $d_2 \cdots d_r$ linear
polynomials
$$ a_{1\,i_2 \cdots i_r}
\quad = \quad 1 \, - \, \sum_{j=2}^{d_1} a_{j \,i_2 \cdots
i_r}. $$
Let $Q$ denote a set of  $d_2 \cdots d_n$ new unknowns
$\, q_{i_2 \cdots i_r i_{r+1} \cdots i_n} \, = \, q_{i_2 \cdots i_n}$,
defining a polynomial ring  $\rr[Q]$.
We introduce the partial factorization map
\begin{equation}
\label{homo}
  \Psi : \rr[D] \rightarrow \rr[A \cup Q]\,, \,\,
p_{i_1 i_2 \cdots i_n}
\, \mapsto \,a_{i_1 \cdots i_r} \cdot  q_{i_2 \cdots i_n}.
\end{equation}
The kernel
of $\,\Psi \,$ is precisely the ideal
$\, I_1 := I_{1 \perpp {\rm nd}(1)|{\rm pa}(1)}$.
Note that
$$ \,q_{i_2 \cdots i_n } \, = \, \Psi(
p_{\, + \,i_2 \cdots i_n} ). $$
Therefore $\Psi$ becomes an epimorphism
if we localize $\rr[D]$ at the product ${\bf p}_1$
of the $\,p_{\, + \,i_2 \cdots i_n} \,$
and we localize $R$ at the product of
the  $q_{i_2 \cdots i_n }$. This implies
that any ideal  $L$ in the polynomial ring $\rr[D]$ satisfies
the identity
\begin{equation}
\label{pushpull}
  \Psi^{-1}(\Psi(L)) \quad = \quad
\bigl((L + I_1) : {\bf p}_1^\infty\bigr).
\end{equation}

Let $G'$ denote the graph obtained from $G$ by removing
the sink $1$ and all edges incident to $1$.
We regard $\,I_{{\rm local}(G')} \,$ as an ideal
in $\rr[Q]$. We modify the set of independence statements
${\rm local}(G)$  by removing
$1$ from the sets ${\rm nd}(i)$ for any $i \geq 2$.
Let $J \subset \rr[D]$ be the ideal corresponding to these modified
independence statements, so that $\Psi(J) = I_{{\rm local}(G')}$.
Note that
$$J + I_1 \,\,\subseteq \,\,
I_{{\rm local}(G)} \,\, \subseteq \,\, I_{{\rm global}(G)} \,\,
\subseteq \,\, \ker (\Phi),$$
so it suffices to show that $(J + I_1) : {\bf p}^\infty = \ker
(\Phi)$. The map $\Phi$ factors as
\begin{equation}
\rr[D] \stackrel{\Psi}{\longrightarrow} \rr[A \cup Q]
\stackrel{\Phi'}{\longrightarrow} \rr[A \cup E'] = \rr[E],
\end{equation}
where $\Phi'$ is the factorization map coming from the graph $G'$,
extended
to be the identity on the variables $A$.
By induction on the number of vertices, we may assume
that Theorem \ref{prob812solved} holds for the smaller graph $G'$,
i.e.,
\begin{equation}
\label{yetanotherkernel}
  {\rm ker}(\Phi') \quad = \quad
(I_{{\rm local}(G')} : {\bf q}_2^\infty) \quad = \quad \Psi(J : {\bf
p}_2^\infty),
\end{equation}
where ${\bf q}_2 = \Psi({\bf p}_2)$ and ${\bf p}_2$ is the product of
the linear forms
$\,p_{++\cdots+u_i\cdots u_n}\,$ with at least two initial $+$'s.
Therefore
\begin{equation}
\ker (\Phi) \quad = \quad \Psi^{-1}(\Psi(J : {\bf p}_2^\infty)).
\end{equation}
Applying (\ref{pushpull}), we get
$\, \ker (\Phi) \, = \, ((J : {\bf p}_2^\infty) + I_1) : {\bf
p_1}^\infty
\, = \, (J + I_1) : {\bf p}^\infty $.
\end{proof}

By following the technique of the proof, we can replace
${\bf  p}_1$ by the product of a much smaller number of
$p_{+u_2\cdots u_n}$.  In  fact,
we need only take the linear forms
$\,p_{+u_2\cdots u_r11\cdots1}$.
Hence, by induction,   ${\bf p}$ can be replaced by
a much smaller product of linear forms.
This observation proved to be crucial for computing
some of the tough primary decompositions in Section 6.

As a corollary  we derive an algebraic proof of the Factorization
Theorem.

\vskip .1cm

\noindent {\sl Proof of Theorem \ref{factheorem}:}
We use induction on the number of nodes
 to show that every point in
$\,V_{\geq}(I_{\local(G)} + \langle p -1 \rangle) \,$
also lies in $\,{\rm image}(\phi_{\geq})$.
Such a point is  a homomorphism
$\, \tau \,: \, \rr[D] \rightarrow  \rr\,$ with the
property that $\tau$ is zero on $I_{\local(G)}$,
and its values on the indeterminates $p_{u_1 \cdots u_n}$
are non-negative and sum to $1$.
The map $\tau$  can be extended  to a  homomorphism
$\,\tau' : \rr[Q \cup A] \rightarrow \rr$ as follows.
We first set
$ \tau'(q_{i_2 \cdots i_n })\, = \, \tau (p_{\, + \,i_2 \cdots i_n} ) $.
If that real number is  positive then we set
$\,\tau'(a_{i_1 \cdots i_r })\, = \,
\tau(p_{i_1 i_2 \cdots i_n} )/\tau (p_{\, + \,i_2 \cdots i_n})$,
and otherwise we set  $\,\tau'(a_{i_1 \cdots i_r })\, = \,  0 $.
Our non-negativity hypothesis implies that
$\,\tau  \, $ coincides with the composition of
$\,\tau' \,$ and $\,\Psi $, i.e., the point $\tau$
is the image of  $\tau'$ under the
induced map $\rr^{A \cup Q} \rightarrow \rr^D$.
The conclusion now follows by induction.
\qed

\vskip .2cm

We close this section by presenting our solution
to \cite[Problem 8.11]{S}.

\begin{prop}
\label{fivenonrad} There exists a Bayesian network $G$
on five binary random variables such
that the local Markov ideal $\,I_{\rm local}(G)\,$
is not radical.
\end{prop}

\begin{proof}
Let $G$ be the complete bipartite
network $\,K_{2,3}\, $ with nodes $\{1,5\}$ and $\{2,3,4\}$
and directed edges $\, (5,2), \, (5,3), \,(5,4) ,
(2,1) , \, (3,1) ,\, (4,1)$. Then
$$ {\rm local}(G) \,\, = \,\,
\bigl\{\, 1 \perpp 5 \,|\, \{2,3,4\} , \,
2 \perpp \{3,4\} \,|\, 5 , \, 3 \perpp \{2,4\} \,|\, 5, \,
4 \perpp \{2,3\} \,| \, 5\,\bigr\} . $$
The polynomial ring $\, \rr[E]\,$ has $32$ indeterminates
$\,p_{11111},  p_{11112}, \ldots, p_{22222}$.
The ideal $\,I_{\rm local}(G)\,$ is minimally generated by
eight binomial quadrics
$$\, p_{1 u_2 u_3 u_4 1}\cdot p_{2 u_2 u_3 u_4 2} \,\, - \,\,
  p_{1 u_2 u_3 u_4 2} \cdot p_{2 u_2 u_3 u_4 1} ,
\qquad u_2,u_3,u_4 \in \{1,2\}, $$
and eighteen non-binomial quadrics
\begin{eqnarray*}
& \!\!\!\! p_{+122u_5} \cdot p_{+221u_5} - p_{+121u_5} \cdot
p_{+222u_5}  , \,
p_{+212u_5} \cdot p_{+221u_5} - p_{+211u_5} \cdot p_{+222u_5}, \\
& \!\!\!\! p_{+112u_5} \cdot p_{+221u_5} - p_{+111u_5} \cdot
p_{+222u_5}   , \,
p_{+122u_5} \cdot p_{+212u_5} - p_{+112u_5} \cdot p_{+222u_5} , \\
& \!\!\!\! p_{+121u_5} \cdot p_{+212u_5} - p_{+111u_5} \cdot
p_{+222u_5}   , \,
p_{+122u_5} \cdot p_{+211u_5} - p_{+111u_5} \cdot p_{+222u_5} ,   \\
& \!\!\!\! p_{+112u_5} \cdot p_{+211u_5} - p_{+111u_5} \cdot
p_{+212u_5}  , \,
p_{+121u_5} \cdot p_{+211u_5} - p_{+111u_5} \cdot p_{+221u_5} ,   \\
& p_{+112u_5} \cdot p_{+121u_5} - p_{+111u_5} \cdot p_{+122u_5}\,,
\qquad \qquad u_5 \in \{1,2\}.
\end{eqnarray*}
These nine equations
(for fixed value of $u_5$) define the Segre embedding
of  $\, \pp^1 \times  \pp^1 \times  \pp^1 \,$ in $\,\pp^7 $,
as in \cite[eqn.~(8.6), page 103]{S}.
Consider the polynomial
$$ f \quad = \quad
p_{+1112} p_{+2222}
(  p_{12221} p_{12212} p_{12122} p_{12111}
  - p_{12112} p_{12121} p_{12211} p_{12222} ). $$
By computing a Gr\"obner basis, it can be checked
that $\, f^2 \,$ lies in $\, I_{{\rm local}(G)}\,$
but $\,f \,$ does not lie in $\, I_{{\rm local}(G)}$.
Hence  $\, I_{{\rm local}(G)}\,$ is not
a radical ideal. The primary decomposition of
this ideal will be described in Example \ref{G138}.
\end{proof}

\section{Networks on Four Random Variables}

In  this section we present
the algebraic classification of all Bayesian networks
on four random variables. In the binary case we have the
following result.

\begin{thm}
\label{binary4class}
The local and global Markov ideals of
all Bayesian networks
on four binary variables are radical.
The hypothesis ``binary'' is essential.
\end{thm}

\begin{table}[!hbt]
\begin{center}
\begin{tabular}{| l | l | l | l | l |}
\hline
Index & Information & Network & Local  & Global
\\ \hline
1 & $1,\, 2,\, 1$ & $\{\}, \{1\}, \{1, 2\}, \{1, 2\}$  & prime &
\\ \hline
2 & $2,\, 4,\, 2$ & $\{\}, \{1\}, \{1\}, \{1, 2, 3\}$  & prime &
\\ \hline
3 & $2,\, 4,\, 2$ & $\{\}, \{1\}, \{1, 2\}, \{1, 3\}$  & prime &
\\ \hline
4 & $3,\, 4,\, 6$ & $\{\}, \{1\}, \{1\}, \{1, 2\}$  & prime &
\\ \hline
5 & $4,\, 6,\, 9$ & $\{\}, \{1\}, \{1\}, \{1\}$  & prime &
\\ \hline
6 & $4,\, 16,\, 4$ & $\{\}, \{\}, \{1, 2\}, \{1, 2, 3\}$  & prime &
\\ \hline
7 & $4,\, 16,\, 4$ & $\{\}, \{1\}, \{1, 2\}, \{2, 3\}$  & prime &
\\ \hline
8 & $4,\, 16,\, 4$ & $\{\}, \{1\}, \{2\}, \{1, 2, 3\}$  & prime &
\\ \hline
9 & $5,\, 32,\, 5$ & $\{\}, \{\}, \{1, 2\}, \{1, 2\}$  & prime &
\\ \hline
10 & $5,\, 32,\, 5$ & $\{\}, \{1\}, \{1, 2\}, \{2\}$  & prime &
\\ \hline
11 & $6,\, 8,\, 10$ & $\{\}, \{1\}, \{1\}, \{2\}$ & radical, 5 comp. &
prime
\\ \hline
12 & $6,\, 16,\, 12$ & $\{\}, \{\}, \{1\}, \{1, 2, 3\}$  & prime &
\\ \hline
13 & $6,\, 16,\, 12$ & $\{\}, \{\}, \{1, 2\}, \{2, 3\}$ & prime &
\\ \hline
14 & $6,\, 16,\, 12$ & $\{\}, \{1\}, \{2\}, \{2, 3\}$  & prime &
\\ \hline
15 & $6,\, 64,\, 6$ & $\{\}, \{1\}, \{1\}, \{2, 3\}$ & radical, 5 comp.
& radical
\\ \hline
16 & $6,\, 64,\, 6$ & $\{\}, \{1\}, \{1, 2\}, \{3\}$ & radical, 9 comp.
& prime
\\ \hline
17 & $6,\, 64,\, 6$ & $\{\}, \{1\}, \{2\}, \{1, 3\}$ & radical, 5 comp.
& radical
\\ \hline
18 & $7,\, 8,\, 14$ & $\{\}, \{1\}, \{2\}, \{3\}$ & radical, 3 comp. &
prime
\\ \hline
19 & $7,\, 8,\, 28$ & $\{\}, \{\}, \{1\}, \{1, 3\}$ & prime &
\\ \hline
20 & $7,\, 24,\, 16$ & $\{\}, \{\}, \{1\}, \{1, 2\}$ & prime &
\\ \hline
21 & $7,\, 32,\, 13$ & $\{\}, \{1\}, \{2\}, \{2\}$ & prime &
\\ \hline
22 & $8,\, 14,\, 31$ & $\{\}, \{\}, \{1\}, \{1\}$ & prime &
\\ \hline
23 & $8,\, 34,\, 20$ & $\{\}, \{\}, \{1\}, \{2, 3\}$ & prime &
\\ \hline
24 & $8,\, 36,\, 18$ & $\{\}, \{\}, \{\}, \{1, 2, 3\}$ & prime &
\\ \hline
25 & $8,\, 36,\, 18$ & $\{\}, \{\}, \{1, 2\}, \{3\}$ & prime &
\\ \hline
26 & $9,\, 20,\, 27$ & $\{\}, \{\}, \{1\}, \{2\}$ & radical, 5 comp. &
prime
\\ \hline
27 & $9,\, 24,\, 34$ & $\{\}, \{\}, \{\}, \{1, 2\}$  & prime &
\\ \hline
28 & $9,\, 24,\, 34$ & $\{\}, \{\}, \{1\}, \{3\}$  & prime &
\\ \hline
29 & $10,\, 20,\, 46$ & $\{\}, \{\}, \{\}, \{1\}$ & prime &
\\ \hline
30 & $11,\, 24,\, 55$ & $\{\}, \{\}, \{\}, \{\}$ & prime &
\\ \hline
\end{tabular}
\end{center}
\caption{All Bayesian Networks on Four Binary Random Variables}
\label{binary4}
\end{table}

Thus the solution
\cite[Problem 8.11]{S} is affirmative for
networks on four binary nodes.
Proposition \ref{fivenonrad}
shows that the hypothesis ``four'' is essential.
Theorem  \ref{binary4class}
is proved by exhaustive computations in   {\tt Macaulay2}.
We summarize the results  in Table \ref{binary4}.
Each row represents one network $G$ on four binary random variables
along with some information about its two ideals
$$ I_{{\rm local}(G)} \,\, \subseteq \,\, I_{\mathrm{global}(G)} \,\,
\subseteq \,\, \rr[p_{1111}, p_{1112}, \ldots, p_{2221}, p_{2222} ]. $$
Here $G$ is represented by the list of sets of children
  $ ({\rm ch}(1),{\rm ch}(2),{\rm ch}(3),{\rm ch}(4)) $.
The information given in the second column corresponds to the
\emph{codimension, degree, and number of minimal generators} of the
ideal $I_{\mathrm{local}(G)}$.
For example, the network in the  fourth row
has four directed edges $(2,1), (3,1), (4,1)$ and $(4,2)$.
Here  $\, I_{{\rm local}(G)} =
I_{{\rm global}(G)} = \ker(\Phi) $. This
prime has codimension $3$, degree $4$ and is
generated by the six $2 \times 2$-minors of the
$2 \times 4$-matrix
$$ \begin{pmatrix}
p_{+111} & p_{+112} & p_{+211} & p_{+212} \\
p_{+121} & p_{+122} & p_{+221} & p_{+222} \\
\end{pmatrix}. $$
Of the $30$ local Markov ideals in Table \ref{binary4}
all but six are prime. The remaining six ideals are
all radical, and the number of their minimal primes is listed.
Hence all local Markov ideals are radical. The last column
corresponds to the ideal $I_{\mathrm{global}(G)}$. This ideal is equal
to the distinguished component for all but two networks, namely
{\tt 15} and {\tt 17}. For these two networks we have
$\,I_{{\rm local}(G)} = I_{\mathrm{global}(G)}$.
This proves the first assertion of
Theorem \ref{binary4class}.

The main point of this section is the second
sentence in Theorem \ref{binary4class}.
Embedded components can
appear when the number of levels increases.
In the next theorem we let  $d_1,d_2,d_3$ and $d_4$
be arbitrary positive integers.

\begin{thm}
\label{binaryfour}
Of the $30$ local Markov ideals on four random variables, $22$
are always prime,  five are not prime but always radical
(numbers {\tt 10},{\tt 11},{\tt 16}, {\tt 18},{\tt 26}
in Table \ref{binary4})
and three are not radical (numbers {\tt 15},{\tt 17},{\tt 21}
in Table \ref{binary4}).
\end{thm}

\begin{proof}
We prove this theorem by an exhaustive case analysis
of all thirty networks.
In most cases, the ideal $I_{{\rm local}(G)}$
can be made binomial by a suitable coordinate
change, just like in the proof of Theorem
\ref{threemodel}. In fact, let us start
with a non-trivial case which is immediately
taken care of by Theorem \ref{threemodel}.

\vskip .1cm \noindent {\sl The network {\tt 16}: }
Here we have
${\rm local}(G)  =
\bigl\{ 1 \perpp 4 \,|\, \{2,3\} ,\,
2 \perpp 4 \, | \, 3 \bigr\}$.
For fixed value of the third node
we get the model $\{ 1 \perpp 4 \,|\,  2
 , \, 4 \perpp 2  \} $
whose ideal was shown to be radical in Theorem
\ref{threemodel}. Hence $I_{{\rm local}(G)}$
is the ideal generated by $d_3$
copies of this radical ideal in disjoint sets of variables.
We conclude that  $I_{{\rm local}(G)}$
is radical and has  $(2^{d_2}-1)^{d_3}$
minimal primes.

\vskip .1cm \noindent {\sl The networks
{\tt 1},
{\tt 2},
{\tt 3},
{\tt 4},
{\tt 6},
{\tt 7},
{\tt 8},
{\tt 12},
{\tt 13},
{\tt 14}: }
In each of these ten cases,
the ideal $I_{{\rm local}(G)}$ is generated by
quadratic polynomials corresponding to
a single conditional independence statement.
This observation implies that
$I_{{\rm local}(G)}$ 
is a prime ideal,
by \cite[Lemma 8.2]{S}.

\vskip .1cm \noindent {\sl The network {\tt 5}: }
Here ${\rm local}(G) $ specifies
the model of complete independence
for the random variables $X_2, X_3 $ and $X_4$.
This means that $ I_{{\rm local}(G)}$
is the ideal of a Segre variety, which
is prime and has a quadratic Gr\"obner basis.

\vskip .1cm \noindent {\sl The networks {\tt 24} and {\tt 25}: }
Each of these two  networks describes the join of $d_{4}$ and $d_{3}$
Segre varieties.
The same reasoning as  in case {\tt 5} applies.

\vskip .2cm \noindent {\sl The network  {\tt 23}: }
Observe that $I_{{\rm local}(G)} =
I_{\mathrm{global}(G)} = I_{1 \perpp \{2,4\} \mid 3} + I_{2 \perpp
  \{1,3\} \mid 4}$. Since $G$ is a directed tree,
Theorem \ref{treethm} implies
that $I_{{\rm global}(G)}$ coincides
with the distinguished prime ideal ${\rm ker}(\Phi)$. Therefore,
$I_{{\rm local}(G)}$ is always prime.

\vskip .1cm \noindent {\sl The networks
{\tt 19}, {\tt 22}, {\tt 27}, {\tt 28}, {\tt 29}, {\tt 30}: }
Each of these six networks has an isolated vertex.
This means that $I_{{\rm local}(G)}$
is the ideal of the Segre embedding of
the product of two smaller varieties
namely, the  projective
space $\pp^{d_i-1}$
corresponding to the isolated vertex $i$ and the
scheme specified by the local
ideal of the remaining network on
three nodes.
The latter ideal is prime and
has a quadratic Gr\"obner basis,
by Proposition \ref{threeallprime},
and hence so is $I_{{\rm local}(G)}$.

\vskip .2cm \noindent {\sl The network  {\tt 20}: }
The ideal $I_{{\rm local}(G)}$ is binomial in
the coordinates $p_{ijkl}$ with $i \in \{+,2,\ldots,d_1\}$.
Generators are
$p_{i_1 j_{2} k l}   p_{i_2 j_{1} k l} - p_{i_1 j_{1} k l}  p_{i_2
  j_{2} k l}$,
$p_{i_1 j_{2} k_{1} l}   p_{i_2 j_{1} k_{2} l}
$ $ - p_{i_1 j_{1} k_{1} l}  p_{i_2 j_{2} k_{2} l}
$, and $p_{+ j_1 k_2 l_{1}}   p_{+ j_2 k_1 l_{2}} -   p_{+ j_1 k_1 l_{1}}
p_{+ j_2 k_2 l_{2}}$.
The S-pairs within each group  reduce to zero by the
Gr\"obner basis property of the $2 \times 2$-minors
of a generic matrix. It can be checked easily that
the crosswise reverse lexicographic S-pairs also reduce to
zero.  We conclude that the given set of
irreducible quadrics is a reverse
lexicographic Gr\"obner basis. In view
of  \cite[Lemma 12.1]{GBCP}, the lowest variable
is not a zero-divisor, and hence by symmetry
none of the variables $p_{ijkl}$ is zero-divisor.
It now follows from equation (\ref{sateqn})
in Theorem \ref{prob812solved} that
$I_{{\rm local}(G)}$ coincides with the
prime ideal  $\,\ker(\Phi)$.

\vskip .2cm \noindent {\sl The network  {\tt 9}: }
The ideal $I_{{\rm local}(G)}$ is
generated by the quadratic polynomials
$p_{i_1 j_2 k l} p_{i_2 j_1 k l} - p_{i_1 j_1 k l} p_{i_2 j_2 k
  l}$,  \
$p_{+ + k_1 l_2} p_{+ + k_2 l_1}  -  p_{+ + k_1 l_1} p_{+ + k_2 l_2}
$. These generators form a  Gr\"obner basis
in the reverse lexicographic order. Indeed, assuming that
$i_1 < i_2$, $j_1 < j_2$, $k_1 < k_2$, $l_1 < l_2$, the leading
terms are
$p_{i_1 j_2 k l} p_{i_2 j_1 k l}$  and
$p_{1 1 k_1 l_2} p_{1 1 k_2 l_1}$.
Hence no leading term from the first group
of quadrics shares a variable with a leading
term from the second group.
Hence the crosswise S-pairs
reduce to zero by  \cite[Prop.~4, \S 2.9]{CLO}.
The S-pairs within each group  also reduce to zero by the
Gr\"obner basis property of the $2 \times 2$-minors
of a generic matrix. Hence the generators
are a Gr\"obner basis. Since the leading
terms are square-free, we see that
the ideal is radical.
An argument similar to the previous case
shows that $I_{{\rm local}(G)}$ is  prime.

\vskip .1cm \noindent {\sl The network  {\tt 18}: }
Here $G$ is a directed chain of length four.
We claim that $ I_{{\rm local}(G)}$ is
the  irredundant intersection of
$2^{d_2} - 1$ primes, and it has
a Gr\"obner basis consisting of
square-free binomials of degree two, three and four.
We give an outline of the proof.
We first turn $I_{{\rm local}(G)} $
into a binomial ideal by taking the coordinates
to be $p_{ijkl}$ with $i \in \{+,2,3,\ldots,d_1\}$.
The minimal primes  are indexed by proper subsets of $[d_2]$.
For each such subset $\sigma$ we introduce the monomial prime
$ M_\sigma =  \langle
  p_{+jkl} \,:\, j \in \sigma,
k \in [d_3], l \in [d_4] \rangle $
and the complementary monomial
$ m_\sigma  =
\prod_{j\in [d_2]\backslash \sigma}
\prod_{k \in [d_3]} \prod_{l \in [d_4]} p_{+jkl}, $
and we define the ideal
$ P_{\sigma}  =  \bigl(
(I_{{\rm local}(G)} \, + \, M_\sigma)\, :\, m_\sigma^\infty
\bigr). $
These ideals are prime, and the union of
their varieties is irredundant and
equals the variety of $I_{{\rm local}(G)}$.
Using Buchberger's S-pair criterion, we
check that the following four types
of square-free binomials are a  Gr\"obner basis:
\begin{itemize}
\item the generators $
p_{i_1 j k_1 l_1} p_{i_2 j k_2 l_2}-
p_{i_1 j k_2 l_2} p_{i_1 j k_2 l_2}$
encoding $ 1 \perpp \{3,4\} \,| \, 2 $,
\item the generators $
p_{+ j_1 k l_1}  p_{+ j_2 k l_2}   -
p_{+ j_1 k l_2}  p_{+ j_2 k l_1}  $
encoding $ 2 \perpp 4 \,|\, 3 $,
\item the cubics $
( p_{+ j_1 k l_1} p_{i j_2 k l_2}  -
  p_{+ j_1 k l_2} p_{i j_2 k l_1}  )
\cdot p_{+ j_2 k_3 l_3}$,
\item the quartics $
(p_{i_1 j_1 k l_1}  p_{i_2 j_2 k l_2} -
    p_{i_1 j_1 k l_2}  p_{i_2 j_2 k l_1} )\cdot
p_{+ j_1 l_3 k_3} \cdot p_{+ j_2 l_4 k_4} $.
\end{itemize}

\vskip .2cm \noindent {\sl The network  {\tt 10}:}
The ideal $I_{{\rm local}(G)}$ is generated by
$p_{i_1 j k l_{2}} p_{i_2 j k l_{1}} - p_{i_1 j k l_{1}} p_{i_2 j k l_{2}}$
and
$\,p_{+ + k_1 l_2} p_{+ + k_2 l_1}  -  p_{+ + k_1 l_1} p_{+ + k_2 l_2}$.
In general, this ideal is not prime, but it is always radical.
If $d_{4} = 2$ then the ideal is always prime, If $d_{4} > 2$, $I_{{\rm local}(G)}$ is the
intersection of the distingushed component and $2^{d_{3}-1}$ prime ideals
indexed by all proper subsets $\sigma \subset [d_{3}]$ as in the
previous network.

\vskip .2cm \noindent {\sl The network  {\tt 11}: }
Here, ${\rm local}(G) =
\bigl\{
1 \perpp 4 \mid \{2,3\}  , \: 2 \perpp 3 \mid 4  , \:
3 \perpp \{2, 4\} \bigr\}$.
The ideal $I_{{\rm local}(G)}$ is binomial in
the coordinates $p_{ijkl}$ with $i \in \{+,2,\ldots,d_1\}$.
It is generated by the binomials $
p_{i_1 j k l_1}   p_{i_2 j k l_2} - p_{i_1 j k l_2}  p_{i_2 j k l_1}
$, $
p_{+ j_1 k_1 l_{1}}   p_{+ j_2 k_2 l_{2}} -   p_{+ j_1 k_2 l_{1}}
p_{+ j_2 k_1 l_{2}}$ encoding the first and third independent statements.
The minimal primes  are indexed by pairs of proper subsets of $[d_2]$
and $[d_{3}]$.
For each such pair of subsets $(\sigma,\tau)$ we introduce the monomial prime
$ M_{(\sigma,\tau)} = \langle
 p_{+jkl} \,:\, j \in \sigma,
k \in \tau, l \in [d_4] \rangle $
and the complementary monomial
$ m_{(\sigma,\tau)} =
\prod_{j\in [d_2]\backslash \sigma}
\prod_{k \in [d_3]\backslash \tau} \prod_{l \in [d_4]} p_{+jkl}, $
and we define the ideal
$ P_{(\sigma,\tau)}  =  \bigl(
(I_{{\rm local}(G)} \, + \, M_{(\sigma,\tau)})\, :\, m_{(\sigma,\tau)}^\infty
\bigr).$
These ideals are prime, and the union of
their varieties equals the variety of $I_{{\rm local}(G)}$.
Moreover, the ideal $I_{{\rm local}(G)}$ is equal to the
intersection of the minimal primes which are indexed by the following pairs:
For each proper $\tau \subset [d_{3}]$ the pair $(\emptyset,\tau)$,
and for each nonempty proper $\sigma \subset [d_{2}]$ the pairs $(\sigma,
\tau)$ where $\tau \subset [d_{3}]$ is any subset of cardinality at
most $d_{3}-2$. In particular, for $d_{2} = d_{3} = 3$, and arbitrary
$d_{1}, d_{4}$, the ideal $I_{{\rm local}(G)}$ has $31$ prime components. For
$d_{2}=2, d_{3}=4$, $I_{{\rm local}(G)}$ has $37$ prime components,
and for $d_{2}=4, d_{3}=2$, $I_{{\rm local}(G)}$ has $17$ prime components.

\vskip .2cm \noindent {\sl The network  {\tt 26}: }
The ideal $I_{{\rm local}(G)}$ is a radical ideal.
The minimal primes  are indexed by all pairs of proper subsets of $[d_3]$
and $[d_{4}]$.
For each such pair $(\sigma,\tau)$ we introduce the
monomial primes
$M_{\sigma} = \langle
 p_{+jkl} \,:\, k \in \sigma, j \in [d_{2}],
l \in [d_{4}] \rangle$, $M_{\tau} = \langle p_{i+kl} \,:\,
l \in \tau, i \in [d_1], k \in [d_{3}]  \rangle $, and
$ M_{(\sigma,\tau)} = M_{\sigma} + M_{\tau}$. Just as before, we introduce
the complementary monomial
$ m_{(\sigma,\tau)}$, and
the ideal $ P_{(\sigma,\tau)} =  \bigl(
(I_{{\rm local}(G)} \, + \, M_{(\sigma,\tau)})\, :\, m_{(\sigma,\tau)}^\infty
\bigr).$
The ideal $I_{{\rm local}(G)}$ is equal to the
intersection of all these prime ideals.

\vskip .2cm \noindent {\sl The network  {\tt 21}: }
Here, ${\rm local}(G) =
\bigl\{
1 \perpp \{3,4\} \,| \, 2  , \,
 3 \perpp 4 \bigr\}$.
The ideal $I_{{\rm local}(G)}$  is generated by the binomials $
p_{i_1 j k_{2} l_2}   p_{i_2 j k_{1} l_1} - p_{i_1 j k_{1} l_1}  p_{i_2 j k_{2} l_2}
$, and the polynomials $
p_{+ + k_1 l_{2}}   p_{+ + k_2 l_{1}} -   p_{+ + k_1 l_{1}}   p_{+ + k_2 l_{2}}$.
This ideal is not radical, in general.
The first counterexample occurs for the case
$ d_1 = d_2 = d_3 = 2$ and $ d_4 = 3 $.
Here $I_{{\rm local}(G)}$ is generated by
$33$ quadratic polynomials in $24$ unknowns.
The degree reverse lexicographic Gr\"obner basis
of this ideal consists of
$123$ polynomials of degree up to $8$.
In this case, $I_{{\rm local}(G)}$ is the intersection of the
distinguished component and the $P$-primary ideal $Q = I_{1 \perpp
  \{3,4\} \,| \, 2} + P^{2}$, where $P$ is the prime ideal generated
by the $12$ linear forms $p_{+ j k l}$.

\vskip .2cm \noindent {\sl The networks  {\tt 15} and {\tt 17}: }
Here, after relabeling network {\tt 17}, ${\rm local}(G) =
\bigl\{ 1 \perpp 4 \,| \, \{2,3\}  , \,
 2 \perpp 3 \,|\, 4 \bigr\}$.
The ideal $I_{{\rm local}(G)}$ is binomial in
the coordinates $p_{ijkl}$ with $i \in \{+,2,\ldots,d_1\}$.
It is generated by the binomials $
p_{i_1 j k l_1}   p_{i_2 j k l_2} - p_{i_1 j k l_2}  p_{i_2 j k l_1}
$, $
p_{+ j_1 k_1 l}   p_{+ j_2 k_2 l} -   p_{+ j_1 k_2 l}   p_{+ j_2 k_1 l}$.
This ideal is not radical, in general.
The first counterexample occurs for the case
$ d_1 = 2$ and $ d_2 = d_3 = d_4 = 3 $.
Here $I_{{\rm local}(G)}$ is generated by
$54$ quadratic binomials in $54$ unknowns.
The reverse lexicographic Gr\"obner basis  consists of
$13,038$ binomials of degree up to $14$.
One of the elements in the Gr\"obner basis is
$$  p_{+111}  p_{+223}  (p_{+331})^2 \cdot
\bigl( p_{2122} p_{2133} p_{2323} p_{2332}
 - p_{2333} p_{2322} p_{2132} p_{2123} \bigr). $$
Removing the square from the third factor, we
obtain a polynomial $f$ of degree $7$ such that
that $f \not\in I$ but $f^2 \in I$.
This proves that $I$ is not radical.
The number of minimal primes of $I_{{\rm local}(G)}$
is equal to $2^{d_2} +  2^{d_3}  - 3 $.
\end{proof}

In the $22$ cases where $I_{\rm local}$ is prime,
it follows from Theorem \ref{prob812solved}
that the global Markov ideal $I_{\rm global}$
is prime as well. Among the remaining cases,
we have  $\,I_{{\rm local}(G)} =  I_{{\rm global}(G)}\,$
for networks $ {\tt 10}, {\tt 15},  {\tt 17}, {\tt 21}$,
and we have  $\,I_{\rm local} \not=  I_{\rm global} \,= \,\ker(\Phi)\,$
for networks
$\,{\tt 11}, {\tt 16}, {\tt 18}, {\tt 26}$.
This discussion implies:

\begin{cor}
\label{binaryfourglobal}
Of the $30$ global Markov ideals on four random variables, $26$
are always prime,  one is not prime but always radical
(number {\tt 10} in Table \ref{binary4})
and three are not radical (numbers {\tt 15},{\tt 17},{\tt 21}
in Table \ref{binary4}).
\end{cor}

It is instructive to examine the
distinguished prime ideal $ P  =  {\rm ker}(\Phi)$
in the last case {\tt 15}, {\tt 17}.
Assume for simplicity
that $d_1 = 2$ but $d_2,d_3$ and $d_4$ are
arbitrary positive integers.
We rename the unknowns 
$ x_{jkl}  =  p_{2jkl}$ and
$ y_{jkl}  =  p_{+jkl}$.
Then we can take $\Phi$ to be  the following
monomial map:
\begin{equation}
\label{nothreeway}
\rr[x_{jkl}, y_{jkl}] \rightarrow
\rr[u_{jk},v_{jl}, w_{kl}]  ,\,
x_{jkl} \mapsto u_{jk} v_{jl} w_{kl}  ,\,
y_{jkl} \mapsto  v_{jl} w_{kl}  ,\,
\end{equation}
 For example, for $d_2=d_3=3$ and $d_4=2$,
the ideal $P = {\rm ker}(\Phi)$ has $361$ minimal generators,
of degrees ranging from two to seven. One generator is
$$
x_{111}  x_{132}  x_{222}  x_{312}  x_{321}  y_{221} y_{331}
 -
x_{112}  x_{131}  x_{221}  x_{311}  x_{322}  y_{232}  y_{321}.
$$
Among the $361$ minimal generators, there are precisely $15$ which
do not contain any variable $y_{ijk}$, namely,
there are nine quartics and six sextics like
$$ x_{112} x_{121} x_{211} x_{232} x_{322} x_{331}
- x_{111} x_{122} x_{212} x_{231} x_{321} x_{332}.$$
These $15$ generators form the Markov basis
for the $3 \times 3 \times 2$-tables in the
no-three-way interaction model.
See \cite[Corollary 14.12]{GBCP} for
a discussion.

The ideal for the no-three-way interaction model
of $d_2 \times d_3 \times d_4$-tables always coincides with
the elimination ideal $ P \cap  \rr[x_{ijk}]$
and, moreover, every generating set of $P$ contains
a generating set for $ P \cap  \rr[x_{ijk}]$.
In view of \cite[Proposition 14.14]{GBCP}, this
shows that the  maximal degree among minimal generators
of $P$ exceeds any bound as $d_2, d_3, d_4$ increases.
In practical terms, it is  hard
to compute these generators even for $d_2 = d_3 = d_4 = 4$.
We refer to the web page
{\tt http://math.berkeley.edu/$\sim$seths/ccachallenge.html}.

\section{Networks on Five Binary Random Variables}

In this section we discuss the global Markov ideals of all 
Bayesian networks on five binary random variables. 
In each case we computed the primary decomposition.
In general, the built-in
primary decomposition algorithms in current computer algebra systems
cannot
compute the primary decompositions of most of these ideals.  In 
the Appendix, we outline some  techniques that allowed us to compute 
these decompositions.  The primary decompositions of the local Markov
ideals
of these networks could also be computed, but they have less regular
structure and are in general more complicated.

There are 301 distinct non-complete networks on five random variables, up to
isomorphism of directed graphs. 
 We have placed descriptions of these networks and
of
the primary decompositions of their global Markov ideals on the website
{\tt http://math.cornell.edu/$\sim$mike/bayes/global5.html}.
In this section, we refer to the graphs as $G_0, G_1, \ldots, G_{300}$,
the indices matching the information on the website.
We summarize our results in a theorem.

\begin{thm}
Of the $301$ global Markov ideals on five binary random variables, 
$220$ are prime, $68$ are radical but not prime, and $13$ are not radical.
\end{thm}

\begin{proof}
The proof is via direct computation with each of these ideals
in  {\tt Macaulay2}. Some of these require little or no computation: if
$G$ is a directed forest, or if there is only one independence statement, 
then the ideal is prime.  Others require substantial computation and some
ingenuity to find the primary decomposition.
Results are posted at the website cited above.

To prove primality, it suffices to compute the ideal quotient of $I =
I_{\gl(G)}$ with respect to a small subset of the
 $\,p_{+++u_r\cdots u_n}$.
Alternatively, one may birationally project $I$ by eliminating
variables, as
in Proposition~\ref{prop:birat}.  In either case, if a zero divisor $x$
is found,
the ideal is not prime.  If some ideal quotient satisfies 
$\, (I : x^2) \neq (I:x)$, then $I$ is not radical.
\end{proof}

The numbers of prime components of the $288$
radical global Markov ideals range from $1$
to $39$. The distribution is given
in the following table:

\begin{table}[!hbt]
\begin{center}
\begin{tabular}{| l | l | l | l | l | l | l | l | l | l |}
\hline
\# of components & 1 & 3 & 5 & 7 & 17 & 25 & 29 & 33 & 39
\\ \hline
\# of ideals & 220 & 8 & 41 & 3 & 9 & 1 & 2 & 3 & 1
\\ \hline
\end{tabular}
\end{center}
\end{table}



\begin{thm}
Conjecture \ref{conj:quadrics} is true for 
Bayesian networks $G$ on five binary random variables.
In each of the  $301$ cases, the
distinguished prime ideal $\ker(\Phi)$ is
generated by homogeneous polynomials of
degree at most eight.
\end{thm}

\begin{proof}
We compute the distinguished component 
from $I_{{\rm global}(G)}$
by saturation, and we check the
result by
using the techniques in the Appendix.
The computation of the distinguished component of the 
$81$ non-prime examples yields that
$64$ of these ideals are generated in degrees $\leq 4$,
twelve are generated in degrees $\leq 6$, and five are generated in
degrees $\leq 8$.
\end{proof}

Theorem~\ref{prob812solved} says that we can decide primality or find
the
distinguished component of
$I_{\gl(G)}$
by inverting each of the $p_{+++u_i\cdots u_n}$.  With some care, it is
possible to
reduce this to a smaller set.  Still, the following  is
unexpected.

\begin{prop}
For all but two networks on five binary random variables, $p_{+1111}$ is
a non-zero divisor on $I = I_{\gl(G)}$ if and only if $I$ is prime.  In
all but these two examples, $I$ is radical if and only if
$\, (I : p_{+1111}^2) = (I : p_{+1111})$.
\end{prop}

\begin{proof}
The networks which do not satisfy the given  property are
$G_{201} = \bigl(\{\}, \{1\}, \{1,2\}, \{1,2\}, \{3,4\}\bigr)$
and $G_{214} = \bigl(\{\}, \{1\}, \{1,2\}, \{3\}, \{1,2,4\} \bigr)$.
After permuting the nodes 4, 5, both the local
and global independence statements of $G_{214}$ are the same as those
for $G_{201}$.  The
global independence statements for $G_{201}$ are
$\bigl\{ \{1,2\} \perpp 5 \,|\, \{3,4\}, \; 3 \perpp 4 \,|\, 5\bigr\}$.
The primary decomposition
for the radical ideal $I = I_{\gl(G_{201})}$ is
   $$I \, = \, \ker(\Phi) \cap
(I + P_{++1\bullet\bullet}) \cap (I + P_{++2\bullet\bullet}) \cap (I + P_{++\bullet 1\bullet}) \cap (I +
P_{++\bullet 2\bullet}),$$
where $\ker(\Phi)$ is the distinguished prime component, 
$$ P_{++1\bullet\bullet} \quad = 
\quad \langle p_{++111}, p_{++112}, p_{++121}, p_{++122} \rangle,$$
 and the other three components are defined in an analogous manner.
Therefore, $p_{+1111}$ is a non-zero divisor modulo $I$.  By examining
all 81 non-prime ideals, we see that all except these two have a minimal
prime containing $p_{+1111}$.  The final statement also follows from
direct
computation.
\end{proof}

We have searched for conditions on the network which would characterize
under what conditions the global Markov ideal is prime, or fails to be prime.
Theorem \ref{treethm} states that if the network
is a directed forest, then the global Markov ideal
is prime. Two possible conditions, the first for primality, 
and the second for
non-primality, are close, but not quite right.  We present them, with
their counterexamples, in the following two propositions.

\begin{prop}
\label{wefoundatree}
There is a unique network $G$ on $5$ binary nodes
whose underlying undirected graph
is a tree, but $I_{\gl(G)}$ is not radical.  Every other network 
whose underlying graph is a tree has prime global
Markov~ideal.
\end{prop}

\begin{proof}
The unique network  is $ \,G_{23} \,=\,
 \bigl( \{\}, \{1\}, \{2\}, \{2\}, \{2\} \bigr)$.
Its local and global Markov independent statements coincide
and are equal to
$$\bigl\{1 \perpp \{3,4,5\} \,|\, 2, \; 3 \perpp \{4,5\}, \;
4 \perpp \{3,5\}, \; 5 \perpp \{3,4\}\bigr\}.$$
Computation using {\tt Macaulay2} reveals
    $$I_{\gl(G_{23})}\quad =
\quad  \ker(\Phi) \, \cap \, (I_{\gl(G_{23})} + 
(P_{+\bullet\bullet\bullet\bullet})^2 ), $$
where $\,P_{+\bullet\bullet\bullet\bullet}\,$ is the ideal
generated by the $16$ linear forms    $\, p_{+u_2u_3u_4u_5}$.
Inspecting the  $81$ non-prime ideals shows
that $G_{23}$ is the only example.
\end{proof}

We say that the network $G$ {\it has an induced
$r$-cycle} if there is an induced subgraph $H$ of $G$ with $r$ vertices
which consists of two disjoint directed paths
which share the same start point and end point.

\begin{prop}
Of the $301$ networks on five  nodes, 
$70$ have an induced $4$-cycle or $5$-cycle.
For exactly two of these, the ideal $I_{\gl(G)}$ is prime.
\end{prop}

\begin{proof}
Once again, this follows by examination of the $301$ cases.
The graphs which have an induced 4-cycle but whose global Markov ideal is
prime are
\begin{eqnarray*}
&  G_{265} = \{\{\}, \{1\}, \{1, 2\}, \{1, 2\}, \{2, 3, 4\}\} \\
\hbox{and} & G_{269} = \{\{\}, \{1\}, \{1, 2\}, \{2, 3\}, \{1, 2, 4\}\}.
\end{eqnarray*}
Removing node 2  results in a $4$-cycle.
The local and global Markov statements are 
all the same up to relabeling:
$\,\bigl\{1 \perpp 5 \,|\, \{2,3,4\}, \; 3 \perpp 4 \,|\, 5\bigr\}$.
\end{proof}

There are four graphs with three induced $4$-cycles, namely, $G_{138}$,
$G_{139}$, $G_{150}$,
$G_{157}$.  The first two graphs give rise to the same (global or local)
independence statements, and similarly for the last two.
The ideal $I_{\gl(G_{138})}$ has the most components of any
of the 301 ideals considered in this section.

\begin{ex} 
\label{G138}
The network
$\, G_{138} = \bigl( \{\}, \{1\}, \{1\}, \{1\}, \{2,3,4\} \bigr)\,$
is isomorphic to the one in Proposition \ref{fivenonrad}.
Its ideal $I_{\gl(G_{138})}$
has $207$ minimal primes, and $37$ embedded primes.  Each
of the
$207$ minimal primary components are prime.  We will
describe the  structure of  these
components.

Let $\, F_{i_1i_2i_3} = \det\begin{pmatrix}
     p_{+i_1i_2i_31} & p_{+i_1i_2i_32} \\
     p_{2i_1i_2i_31} & p_{2i_1i_2i_32}
     \end{pmatrix}$.
Let $J_i$ be the ideal generated by the 
$2 \times 2$ minors located in the first two rows
or columns of the matrix
$$\begin{pmatrix}
p_{+111i} & p_{+112i} & p_{+211i} & p_{+212i} \\
p_{+121i} & p_{+122i} & p_{+221i} & p_{+222i} \\
p_{+211i} & p_{+212i} & * & * \\
p_{+221i} & p_{+222i} & * & *
\end{pmatrix}. $$
We have    $$ I \quad := \quad I_{\gl(G_{138})}
\quad = \quad  J_1 \,+\, J_2 \,+\, 
\langle F_{111}, F_{112}, \ldots, F_{222} \rangle .$$
Each $J_i$ is minimally generated by 9 quadrics, so that $I$
is minimally generated by 26 quadrics.
Each $J_i$ is prime of codimension 4, and so $J_1 + J_2$ is prime of
codimension 8.  Since there are only 8 more quadrics, Krull's principal
ideal theorem tells us that all minimal primes have codimension
at most 16, which is also the codimension of the distinguished
component. Note that $I$ is  a binomial ideal
in the unknowns
$\,p_{+u_2 u_3 u_4 u_5}$ and $p_{2 u_2 u_3 u_4 u_5}$.

\begin{table}[!hbt]
\begin{center}
\begin{tabular}{| l | l | l | l |}
\hline
\# primes & codim & degree & faces \\ \hline
6 & 14 & 48 & $(f,f)$, $f$ a facet \\ \hline
12 & 14 & 4 & $(e,e)$, $e$ an edge \\ \hline
24 & 16 & 15 & $(f_1,f_2)$, $f_1 \cap f_2$ is an edge \\ \hline
48 & 16 & 4 & $(f,e)$, $f \cap e$ is a point \\ \hline
12 & 16 & 1 & $(e_1,e_2)$, 2 antipodal edges \\ \hline
48 & 16 & 1 & $(e_1,e_2)$, 2 non-parallel disjoint edges \\ \hline
48 & 16 & 1 & $(e,p)$, point $p$ on the edge antipodal
to $e$ \\ \hline
8 & 16 & 1 & $(p_1,p_2)$, antipodal points \\ \hline
1 & 16 & 2316 & distinguished component \\ \hline
\end{tabular}
\end{center}
\caption{All $207$ minimal primes of the ideal $I_{\gl(G_{138})}$}
\label{ideal138}
\end{table}

Let $\Delta$ be the unit cube, with vertices $(1,1,1), (1,1,2), \ldots,
(2,2,2)$.  If $\sigma \subset \Delta$ is a face, define
$P_{\sigma, i}$ to be the monomial prime generated by
    $\{ p_{+vi} \mid v \not\in \sigma\},$ for $i \in \{1,2\}$.
If $P$ is a minimal prime of $I$, which is not the
distinguished
component, then $P$ must contain some $p_{+v_1v_2v_31}$, and also
contain
some $p_{+u_1u_2u_32}$.  Therefore, there are faces $\sigma_1$
and $\sigma_2$ of $\Delta$ such that 
$P$ contains $P_{\sigma_1,1} + P_{\sigma_2,2}$, and does not
contain any other elements $p_{+vi}$.  Let $m_{\sigma_1\sigma_2}$ be the
product of all of the $p_{+vi}$  such that $v \in \sigma_i$
for $i=1,2$. It turns out that
every minimal prime ideal of $I$ has the form
    $$ P_{\sigma_1,\sigma_2} \quad := \quad \bigl(
 (I + P_{\sigma_1,1} + P_{\sigma_2,2}) :
    m_{\sigma_1\sigma_2}^\infty \bigr) $$
for some pair $\sigma_1, \sigma_2$ of
proper faces  of the cube  $ \Delta$.
However, not all pairs
of faces correspond to
minimal primes.  There are 27 proper faces of the cube, and so there
    are $27^2 = 729$ possible minimal primes.  Only 206 of these occur.
The list of minimal primes is given in Table~\ref{ideal138}. 
\qed
\end{ex}

Bayesian networks give rise to very
interesting  (new and old) constructions in
algebraic geometry.  In the next section, we shall
encounter  secant varieties.  Here, we offer
 a generalization of Example \ref{G138}
to arbitrary toric varieties. 
Let $\,I_A \subset \rr[z_1, \ldots, z_n]$ be any
{\em toric ideal}, specified as in \cite{GBCP}
by a  point configuration $A = \{ a_1, \ldots, a_n \} \subset \zz^d $.
Let $\Delta$ be the convex hull of $A$ in $\rr^d$.
We define the {\em double join} of the toric ideal $I_A$ to be
the new ideal
$$ I_A(x) + I_A(y) + \langle F_1, \ldots, F_n\rangle \, \subset\,
 \rr[x_1, \ldots, x_n, y_1, \ldots,y_n, a_1, \ldots
a_n, b_1, \ldots,
    b_n]$$
where $F_i = \det \begin{pmatrix}
x_i & a_i \\
y_i & b_i
\end{pmatrix}$,
and $I_A(x)$ and $I_A(y)$ are generated by
copies of $I_A$ in 
$ \rr[x_1, \ldots, x_n]$ and $\rr[ y_1, \ldots,y_n]$
respectively.
The ideal $I$ in Example \ref{G138}
is the double join of the Segre variety
$\pp^1 \times \pp^1 \times \pp^1 \subset \pp^7$,
which is the toric variety whose polytope $\Delta$ is the $3$-cube.
In general, the minimal primes of the double join
of $I_A$ are indexed by pairs of faces of
the polytope $\Delta$.
We believe that this construction deserves the
attention of algebraic geometers.

\section{Hidden Variables and Higher Secant Varieties}

Let $G$ be a Bayesian network on $n$ discrete random variables
and let $P_G = {\rm ker}(\Phi)$ be its  homogeneous
prime ideal in the polynomial ring $\rr[D]$,
whose indeterminates $\, p_{i_1 i_2 \cdots i_n}\,$ represent
probabilities of events  $\,(i_1,i_2,\ldots,i_n)  \in D $.
We now consider the situation when some
of the random variables are hidden. After relabeling
we may assume that the variables corresponding
to the nodes $r+1,\ldots,n$ are hidden, while the
random variables corresponding to the nodes
$1,\ldots,r$ are observed. Thus the \emph{observable probabilities} are
$$ p_{i_1 i_2 \cdots i_r \, + +\cdots +}
\quad = \quad
\sum_{j_{r+1} \in [d_{r+1}]}
\sum_{j_{r+2} \in [d_{r+2}]} \cdots
\sum_{j_{n} \in [d_{n}]}
  p_{i_1 i_2  \cdots i_r \, j_{r+1} j_{r+2} \cdots j_n} .$$
We write  $D' = [d_1] \times \cdots \times [d_r]$
and $\rr[D']$ for the polynomial subring
of $\rr[D]$ generated by the observable
probabilities $\,p_{i_1 i_2 \cdots i_r \, + + \cdots +}$.
Let  $\,\pi : \rr^D \rightarrow \rr^{D'}\,$ denote
the canonical linear epimorphism induced by the
inclusion of $\rr[D']$ in $\rr[D]$.
We are interested in the following inclusions of semi-algebraic sets:
\begin{equation}
\label{inclusions}
  \pi(V_{\geq 0}(P_G)) \,\, \subset \,\,
\pi(V(P_G))_{\geq 0} \,\, \subset \,\,
\pi(V(P_G)) \,\, \subset \,\,
\overline{\pi(V(P_G))} \,\,\subset\,\,
\rr^{D'}.
\end{equation}
These inclusions are generally all strict.
In particular, the space $\,\pi(V_{\geq 0}(P_G))$ which
consists of all \emph{observable probability distributions}
is often much smaller than the space
$\,\pi(V(P_G))_{\geq 0} $ which consists of
probability distributions on $D'$ which
would be observable if  non-negative or complex
numbers were allowed for the hidden parameters.
However, they have the same Zariski closure:

\begin{prop}
The set of all polynomial functions which vanish on the space
  $\, \pi(V_{\geq 0}(P_G)) \,$ of
observable probability distributions
is the prime ideal
\begin{equation}
\label{hiddenideal}
  Q_G \quad = \quad  P_G \,\cap \, \rr[D'].
\end{equation}
\end{prop}

\begin{proof}
The elimination ideal $Q_G \subset \rr[D']$ is prime because
$P_G \subset \rr[D]$ was a prime ideal. By the
Closure Theorem of Elimination Theory \cite[Theorem 3, \S 3.2]{CLO},
the ideal $Q_G$ is the vanishing ideal of the image
$\,\pi(V(P_G)) $. Since
$\,V_{\geq 0}(P_G) \,$ is Zariski dense  in $\,V(P_G) $,
by the Factorization Theorem  \ref{factheorem},
 and $\pi$ is a linear map, it follows that
$\,\pi(V_{\geq 0}(P_G)) \,$ is Zariski dense  in $\,\pi(V(P_G)) $.
\end{proof}

We wish to demonstrate how computational algebraic geometry can
be used to study  hidden random variables in
Bayesian networks. To this end we apply the concepts
introduced above to a standard example
from the statistics literature
\cite{GHKM}, \cite{SS1}, \cite{SS2}.
We fix the network $G$ which has
  $n+1$ random variables $F_1,\ldots,F_n, H$ and
$n$ directed edges $(H,F_i)$, $ i = 1,2,\ldots,n$.
This is the \emph{naive Bayes model}. The variable
$H$ is the hidden variable, and its levels
$\,1,2,\ldots, d_{n+1} =: r\,$ are called the \emph{classes}.
The observed random variables $F_1,\ldots,F_n$ are
  the \emph{features} of the model.
In this example, the prime ideal $P_G$ coincides with the
local ideal $I_{{\rm local}(G)}$ which is specified by
requiring that, for each fixed class,
  the features are completely independent:
$$ F_1 \perpp F_2 \perpp \cdots \perpp F_n \,|\, H . $$
This ideal is obtained as the kernel of the map
$\,p_{i_1 i_2 \cdots i_n k} \mapsto  x_{i_1} y_{i_2} \cdots z_{i_n} $,
one copy for each fixed class $k$, and then
adding up these $r$ prime ideals.
Equivalently, $P_G$ is the  ideal of the join of
$r $ copies of the \emph{Segre variety}
\begin{equation}
\label{segre}
X_{d_1,d_2,\ldots,d_n} \quad := \quad
\pp^{d_1-1} \times
\pp^{d_2-1} \times \cdots \times
\pp^{d_n-1} \quad \subset \quad
\pp^{d_1 d_2 \cdots d_n - 1}.
\end{equation}
The points on  $\,X_{d_1,d_2,\ldots,d_n} \,$
represent tensors of rank $ \leq 1$.
Our linear map $\pi$ takes an $r$-tuple
of tensors of rank $ \leq 1$ and it computes
their sum, which is a tensor of rank $\leq r$.
The closure of the image of $\pi$ is
what is called a \emph{higher secant variety}
in the language of algebraic geometry
\cite[Example 11.30]{Harris}.

\begin{cor}
The naive Bayes model with $r$ classes
and $n$ features corresponds to the
$r$-th secant variety
of a Segre product of $n$ projective spaces:
\begin{equation}
\label{segresecant}
\overline{\pi(V(P_G))} \quad  = \quad
{\rm Sec}^{r}(X_{d_1,d_2,\ldots,d_n})
\end{equation}
\end{cor}

The case $n = 2$ of two features is a staple of
classical projective geometry. In that
special case, the image of $\pi$ is closed, and
$\, \pi(V(P_G)) \,= \,{\rm Sec}^{r}(X_{d_1,d_2})\,$ consists of
all real  $d_1 \times d_2$-matrices of
rank at most $r$. This variety has codimension
$\, (d_1-r)(d_2-r)$, provided $r \leq {\rm min}(d_1,d_2)$.
Its ideal $Q_G$ is generated by the $(r+1)\times (r+1)$-minors
of the $d_1 \times d_2$ matrix $(p_{ij+})$. The dimension
formula of Settimi and Smith  \cite[Theorem 1]{SS2}
follows  immediately. For instance,
in the case of two ternary features
$(d_1 = d_2 = 3, r = 2)$, discussed
in different guises in
\cite[\S 4.2]{SS2} and  \cite[Example 11.26]{Harris},
the \emph{observable space} is the
cubic hypersurface defined by the $3 \times 3$-determinant
${\rm det}(p_{ij+})$.

The leftmost inclusion
in (\ref{inclusions}) leads to difficult open problems
even for $n = 2$ features.  Here,
$\,\pi(V(P_G))_{\geq 0} \,$ is the
set of all non-negative $d_1 \times d_2$-matrices
of rank at most $r$, while
$\, \pi(V_{\geq 0}(P_G)) \,$ is the subset consisting
of all matrices of \emph{non-negative rank} at most $r$.
Their difference consists of non-negative matrices of rank $\leq r$
which cannot be written as the sum of $r$ non-negative
matrices of rank $1$. In spite of recent progress
by Barradas and Solis \cite{BS}, there is still no
practical algorithm for computing the
non-negative rank of a $d_1 \times d_2$-matrix.
Things get even harder for $n \geq 3$,
when testing membership in $\, \pi(V_{\geq 0}(P_G)) \,$ means
computing \emph{non-negative tensor rank}.

We next discuss what is known about the case
of $n \geq 3$ features.  The
\emph{expected dimension} of the
secant variety (\ref{segresecant}) is
\begin{equation}
\label{expdim}
  r \cdot (d_1 + d_2 + \cdots + d_n - n + 1) \,\, - \,\,1 .
\end{equation}
This number is always an upper bound, and it
is an interesting problem, studied in
the statistics literature in \cite{GHKM}, to
characterize those cases $(d_1,\ldots,d_n; r)$ when
the dimension is less than the expected dimension.
We note that the results on  dimension in \cite{GHKM}
are all special cases of results by
Catalisano, Geramita and Gimigliano \cite{CGG}, and
the results on singularities in \cite{GHKM} follow from the
geometric fact that the $r$-th secant variety
of any  projective variety is always singular
along the $(r-1)$-st secant variety.
The statistical problem of
\emph{identifiability}, addressed in \cite{SS1},
is related to   the beautiful work
of Strassen \cite{Str} on tensor rank, notably his
Theorem 2.7 on \emph{optimal computations}.

\begin{table}
\begin{center}
$$
\begin{matrix}
dim(X) & dim( {\rm Sec}^{2}(X)) &
\prod_{i=1} d_i & (d_1,\ldots,d_n) & degree & cubics \\
  4  & 9 & 12   & (2,2,3) &   6       &  4 \\
  4  & 9 & 16  & (2,2,2,2) &   64       &  32 \\
  5 & 11 & 16 & (2,2,4) & 20 & 16 \\
  5 & 11 & 18 & (2,3,3) & 57 & 36 \\
  5 & 11 & 24 & (2,2,2,3) & 526 & 184 \\
  5 & 11 & 32 & (2,2,2,2,2) & 3256 & 768 \\
  6 & 13 & 20 & (2,2,5) & 50 & 40 \\
  6 & 13 & 24 & (2,3,4) & 276 & 120 \\
  6 & 13 & 27 & (3,3,3) & 783 & 222 \\
  6 & 13 & 32 & (2,2,2,4) & 2388 & 544 \\
  6 & 13 & 36 & (2,2,3,3) & 6144 & 932
\end{matrix} $$
\end{center}
\caption{The
prime ideal defining the secant lines to the Segre variety (\ref{segre})}
\label{secantlines}
\end{table}

In Table \ref{secantlines} we display
the range of straightforward
{\tt Macaulay2} computations when $\, dim(X) \, = \,
  d_1+\cdots+d_n- 1\,$ is small.
First consider the case  of two classes ($r=2$),
which corresponds to secant lines on $X =
\pp^{d_1-1} \times \cdots \times  \pp^{d_n-1} $.
In each of these cases, the ideal $Q_G$ is generated
by cubic polynomials, and each of these cubic generators
is the determinant of a two-dimensional matrix
obtained by flattening the tensor $\,(p_{i_1 i_2 \cdots i_n})$.
The column labeled ``cubics'' lists the number of minimal
generators. For example, in the case $(d_1=d_2=d_3=3)$, we can
flatten $(p_{ijk})$ in three possible ways to a
$3 \times 9$-matrix, and these have
$\, 3 \cdot \binom{9}{3}\, \,= \, 252 \,$
maximal subdeterminants. The vector space
spanned by these subdeterminants has dimension
$222$, the listed number of minimal generators.
The column ``degree'' lists the degree of
the projective variety $\,{\rm Sec}^2(X)$, which
is $783$ in the previous example.
These computational results  in Table \ref{secantlines}
lead us to make the following conjecture:

\begin{conj}
\label{genbycubics}
The prime ideal $Q_G$ of any naive Bayes model $G$
with $r = 2$ classes is generated by
the $3 \times 3$-subdeterminants of any two-dimensional
table obtained by flattening
the $n$-dimensional table $(p_{i_1 i_2 \cdots i_n})$.
\end{conj}

It was proved by Catalisano, Geramita and Gimigliano
that  the variety $\,{\rm Sec}^2(X)$
always has the expected dimension  (\ref{expdim})
when $r = 2$. A well-known example (see \cite[page 221]{Good})
when the dimension is less
than expected occurs for four classes and
three binary features $\,(r = 3, n = 4,
d_1 = d_2 = d_3 = d_4 = 2)$.
Here (\ref{expdim}) evaluates to  $14$, but
$\,{\rm dim}({\rm Sec}^3(X)) \, = \, 13 \,$ for
$\, X = \pp^1 \times \pp^1 \times \pp^1  \times \pp^1$.
The corresponding ideal $Q_G$ is a complete intersection
generated by any two of the three
$4 \times 4$-determinants obtained by
flattening the $2 \times 2 \times 2 \times 2 $-table
$(p_{ijkl})$. The third is a signed sum of the other two.

The problem of identifying explicit generators of $Q_G$ is
much more difficult when $r \geq 3$, i.e.,
when the hidden variable has three or more levels.
We present the complete solution for the
case of three ternary features. Here
$(p_{ijk})$ is an indeterminate $3 \times 3 \times 3$-tensor
which we wish to write as a sum of $r$ rank one tensors.
The following solution is derived from a result of
Strassen \cite[Theorem 4.6]{Str}. Let
$\,A = (p_{ij1})\,$, $\,B = (p_{ij2})\,$ and $\,C = (p_{ij3})\,$ be
three $3 \times 3$-matrices obtained
by taking slices of the $3 \times 3 \times 3$-table $(p_{ijk})$.

\begin{prop}
Let $Q_G$ be the ideal of $\,{\rm Sec}^{r}(\pp^2 \times
\pp^2 \times \pp^2)$,
the naive Bayes model with $n=3$ ternary features
with $r$ classes.
If $r=2$ then  $Q_G$ is generated by the cubics
described in Conjecture \ref{genbycubics}.
If $r=3$ then $Q_G$ is generated by the quartic
entries of the various $3  \times 3$-matrices of the form
$\, A \cdot {\rm adj}(B) \cdot C
- C \cdot {\rm adj}(B) \cdot A $.
If $r=4$ then $Q_G$ is the principal ideal
generated by the following homogeneous polynomial
of degree $9$ with   $9,216$ terms:
$$ {\rm det}(B)^2 \cdot {\rm det}\bigl( A \cdot B^{-1} \cdot C
-  C \cdot B^{-1} \cdot B \bigr). $$
If $r \geq 5$ then $Q_G$ is the zero ideal.
\end{prop}

\bigskip
\bigskip

\subsection*{Appendix:  Techniques for  Primary Decomposition}

The ideals in this paper present a challenge for present day computer
algebra
systems.  Their large number of variables (e.g. 32
in Section 6), combined
with the
sometimes long polynomials which arise are difficult to handle with
built-in
primary decomposition algorithms.  Even the standard implementations of
factorization of multivariate polynomials have difficulty with some of
the
long polynomials.  This is only a problem with current implementations,
which
are generally not optimized for large numbers of variables.

For the computations performed in Sections 5 and 6,
it was necessary to write
special code (in {\tt Macaulay2}) in order to compute the components
and primary decompositions of these ideals.
We also have some code  in {\tt Macaulay2}
or {\tt Singular}
for generating the ideals $I_{{\rm local}(G)}$ or $I_{{\rm global}(G)}$
from the graph $G$ and the integers
$d_1,d_2,\ldots,d_n$.
In this appendix we indicate
some techniques and tricks that were
used to compute with these ideals.

The first modification which simplifies the problems dramatically is to
change coordinates so that the indeterminates are $p_{2u_2\cdots u_n}$
and
$p_{+u_2\cdots u_n}$, instead of $p_{u_1\cdots u_n}$.  This change
of variables sometimes takes a Markov ideal into a binomial ideal,
which is
generally much simpler to compute with.
Computing  any one Gr\"obner basis, ideal quotient, or
intersection of our ideals is not too difficult.
Therefore, our algorithms make use of these operations.
All ideals examined in this project have the
property that every component is rational.
The distinguished component $\, \ker(\Phi) \,$ is more
complicated than any of the other components, in terms of the number of
generators and their degrees, and it cannot be computed
by implicitization.

The first problem is to decide whether an ideal is prime
(i.e.~whether it equals the unknown ideal $\, \ker(\Phi) $).  There are
several known methods  for deciding primality
(see  \cite{DGP} for a nice exposition).
 The standard method is to reduce to a
zero-dimensional problem.  This entails either a generic
change of coordinates, or factorization over extension fields.
We found that the current implementations of these methods fail for the
majority of the $301$ examples  in Section 6.
The technique that did work for us is to search for
birational projections.
This either produces a zero divisor, or a proof that the
ideal is prime.  It can sometimes be used to count the
components  (both
minimal and embedded), without actually producing the components.

The following result is proved by localizing
with respect to powers of $g$.
This defines a {\em birational projection}
$(x_1,x_2,\ldots,x_n) \mapsto (x_2,\ldots,x_n)$ for $J$.

\begin{prop}~\label{prop:birat}
Let $J \subset  \rr[x_1, \ldots, x_n]$ be an ideal, containing a  polynomial
$f = gx_1 + h$, with $g,h$ not involving $x_1$, and $g$ a non-zero  divisor
modulo $J$.  Let $\, J_1 = J \cap \rr[x_2, \ldots, x_n]\,$ be the
elimination ideal. Then

(a) $J \, = \bigl( \langle J_1, gx_1+h \rangle  : g^\infty \bigr)$,

(b) $J\,$ is prime if and only if $\,J_1\,$ is prime.

(c) $J\,$ is primary if and only if $\,J_1\,$ is primary.

(d) Any irredundant primary decomposition  of $\,J_1\,$
lifts to an irredundant primary decomposition of $J$.
\end{prop}

Our algorithm to check primality starts by searching for variables which
occur linearly, checking that its lead coefficient is not a zero divisor
and then eliminating that variable as in Proposition \ref{prop:birat}.
In almost all of the Markov ideals that we have studied,
iterative use of  this technique proves
or disproves primality.  A priori, one might not be able to find a
birational projection at all, but  this never happened for
any of our examples.

The second problem is to compute the minimal primes or the primary
decomposition. Finding the minimal
primes is the first step in computing a primary  decomposition,
using the technique of \cite{SY}, which is implemented in several
computer algebra systems, including {\tt Macaulay2}.
Here, we have not found a single method  that
always works best. One method that worked in most cases is based on
splitting the ideal into two parts.
Given an ideal $I$, if there is an element $f$ of its Gr\"obner basis
which
factors as $f = f_1f_2$, then
   $$\sqrt{I} \quad = \sqrt{\langle I,f_1 \rangle} \,\,
\cap \,\, \sqrt{\langle I,f_2 \rangle : f_1^\infty}.$$
We keep a list of ideals whose intersection has the same radical as
$I$.  We
process this list of ideals by ascending order on its codimension.  For
each ideal,
we keep a list of the elements that we have inverted by so far (e.g.
$f_1$ in the
ideal $\, \bigl(\langle
I,f_2 \rangle : f_1^\infty \bigr)$)
 and saturate at each step with these
elements.

If there is no element which factors, then we search for a variable to
birationally
project away from, as in Proposition \ref{prop:birat}.  If its lead
coefficient $g$ is a zero
divisor, use this element to split the ideal via
   $$\sqrt{I} \quad = \quad \sqrt{I : g}\,\, \cap\,
\,\sqrt{ \langle I, g \rangle}.$$
As we go, we
 only process ideals which do not contain the intersection of
all known components computed so far.

If we cannot find any birational projection
or reducible polynomial, then
we have no choice but to decompose the ideal  using the
built-in routines, which are based on characteristic sets.  However, in
none
of the examples of this paper was this final step reached.
This method works in a reasonable amount of time
for all but about $10$ to $15$ of the
$301$ ideals in Section 6.

\bigskip
\bigskip
\bigskip
\bigskip

\noindent Luis David Garcia,  Virginia Bioinformatics Institute
at Virginia Tech,
\break Blacksburg, VA 24061, USA,
{\tt lgarcia@vt.edu}.

\bigskip

\noindent Michael Stillman, Department of Mathematics,
Cornell University, Ithaca, NY 14853-4201, USA,
{\tt mike@math.cornell.edu}.

\bigskip

\noindent Bernd Sturmfels, Department of Mathematics,
University of California, \break Berkeley, CA 94720-3840, USA,
{\tt bernd@math.berkeley.edu}.

\end{document}